%Paper: q-alg/9504011
%From: Vitaly Tarasov <tarasov@kusm.kyoto-u.ac.jp>
%Date: Thu, 20 Apr 95 12:25:51 JST

%%%%%%%%%%%%%%%%%%%%%%%%%%%%%%%%%%%%%%%%%%%%%%%%%%%%%%%%%%%%%%%%%%%%%%%%%
%									%
%   Bases of Bethe Vectors and						%
%   Difference Equations with Regular Singular Points			%
%									%
%   Vitaly Tarasov and Alexander Varchenko				%
%									%
%   Kyoto-Math 95-04	( amstex.tex (ver. 2.1) is required )		%
%									%
%%%%%%%%%%%%%%%%%%%%%%%%%%%%%%%%%%%%%%%%%%%%%%%%%%%%%%%%%%%%%%%%%%%%%%%%%

%\mag 1200
\input amstex

\overfullrule 0pt

\hsize 6.25truein
\vsize 9.63truein

\expandafter\ifx\csname langl.def\endcsname\relax \else\endinput\fi
\expandafter\edef\csname langl.def\endcsname{%
 \catcode`\noexpand\@=\the\catcode`\@\space}
\catcode`\@=11

\let\@ft@\expandafter

\newif\ifMag
\ifnum\mag>1000 \Magtrue\fi

\mathsurround 1.6\p@
\font\Bbf=cmbx12 
=\ifMag cmr8\else cmr9\fi

\def\iitem#1{\par\hangindent1.5\parindent
 \hglue-.5\parindent\textindent{\rm#1}}
\def\Item#1{\par\setbox\z@\hbox{\rm #1\enspace}\hangindent\wd\z@
 \hglue-2\parindent\kern\wd\z@\textindent{\rm#1}}

\def\center{\par\begingroup\leftskip\z@ plus \hsize \rightskip\leftskip
 \parindent\z@\parfillskip\z@skip \def\\{\break}}
\def\endcenter{\endgraf\endgroup}

\def\Abstract{\begingroup\narrower\nt{\bf Abstract.}\enspace\ignore}
\def\endAbs{\endgraf\endgroup}

\let\b@gr@@\begingroup \let\B@gr@@\begingroup
\def\b@gr@{\b@gr@@\let\b@gr@@\undefined}
\def\B@gr@{\B@gr@@\let\B@gr@@\undefined}

\def\r@st@re#1{\let#1\s@v@} \def\s@v@d@f{\let\s@v@}

\def\p@sk@p#1#2{\par\skip@#2\relax\ifdim\lastskip<\skip@\relax\removelastskip
 \ifnum#1=\z@\else\penalty#1\relax\fi\vskip\skip@
 \else\ifnum#1=\z@\else\penalty#1\relax\fi\fi}
\def\sk@@p#1{\par\skip@#1\relax\ifdim\lastskip<\skip@\relax\removelastskip
 \vskip\skip@\fi}

\def\s@ct#1#2{\p@sk@p{-200}{.8\bls}\vtop{\bf\setbox\z@\hbox{#1}%
 \parindent\wd\z@\ifdim\wd\z@>\z@\adv\parindent.5em\fi
 \hang\textindent{#1}#2\strut}%
 \nointerlineskip\nobreak\vtop{\strut}\nobreak\vskip-.6\bls\nobreak}

\def\gadv{\global\adv} \def\gad#1{\gadv#1 1} 
\def\l@b@l#1#2{\def\n@@{\csname #2no\endcsname}%
 \if *#1\gad\n@@ \@ft@\xdef\csname @#1@#2@\endcsname{\the\Sno.\the\n@@}%
 \else\@ft@\ifx\csname @#1@#2@\endcsname\relax\gad\n@@
 \@ft@\xdef\csname @#1@#2@\endcsname{\the\Sno.\the\n@@}\fi\fi}
\def\l@bel#1#2{\l@b@l{#1}{#2}\?#1@#2?}
\def\?#1?{\csname @#1@\endcsname}
\def\[#1]{{\count@\z@%
 \setbox\z@\hbox{\?#1@L?}\ifdim\wd\z@>\z@\adv\count@1\fi
 \setbox\z@\hbox{\?#1@F?}\ifdim\wd\z@>\z@\adv\count@1\fi
 \setbox\z@\hbox{\?#1?}\ifdim\wd\z@>\z@\adv\count@1\fi\relax
 \ifnum\count@=\z@{\bf ???}\else\ifnum\count@=1\let\@@ps\relax\let\@@@\relax
 \else\let\@@ps\underbar\def\@@@{{\rm<}}\fi
 \@@ps{\?#1?\@@@\?#1@L?\@@@\?#1@F?}\fi}}
\def\(#1){{\rm(\[#1])}}

\def\dff{\@ft@\d@f} \def\d@f{\@ft@\def}
\def\edff{\@ft@\ed@f} \def\ed@f{\@ft@\edef}

\def\pr@cl#1{\r@st@re\pr@c@\pr@c@{#1}\global\let\pr@c@\relax}

\newcount\Sno \newcount\Lno \newcount\Fno
\def\Section#1{\gad\Sno\Fno\z@\Lno\z@\s@ct{\the\Sno.}{#1}}
\def\section#1{\gad\Sno\Fno\z@\Lno\z@\s@ct\empty{#1}}
\let\Sect\Section \let\sect\section
\def\l@F#1{\l@bel{#1}F} \def\l@L#1{\l@bel{#1}L} \def\<#1>{\l@b@l{#1}F}
\def\Tag#1{\tag\l@F{#1}} \def\Tagg#1{\tag"\rlap{\rm(\l@F{#1})}"}
\def\Th#1{\pr@cl{(\l@L{#1}) Theorem}\ignore}
\def\Lm#1{\pr@cl{(\l@L{#1}) Lemma}\ignore}
\def\Cr#1{\pr@cl{(\l@L{#1}) Corollary}\ignore}
\def\Df#1{\pr@cl{(\l@L{#1}) Definition}\ignore}
\def\Cj#1{\pr@cl{(\l@L{#1}) Conjecture}\ignore}
\def\Rem{\demo{\sl Remark}} 
\def\Pf#1.{\demo{Proof #1}} \def\epf{\qed\enddemo}
\def\qed{\null\nobreak\hfill\nobreak{\=$\,\square$}}
\def\back#1 {\strut\kern-.33em #1\enspace\ignore} % !!! a space after #1 !!!
\def\Text#1{\crcr\noalign{\alb\vsk>\normalbaselines\vsk->\vbox{\nt #1\strut}%
 \nobreak\nointerlineskip\vbox{\strut}\nobreak\vsk->\nobreak}}

\ifx\plainfootnote\undefined \let\plainfootnote\footnote \fi

\let\s@v@\proclaim \let\proclaim\relax
\def\r@R@fs#1{\let#1\s@R@fs} \let\s@R@fs\Refs \let\Refs\relax
\def\r@endd@#1{\let#1\s@endd@} \let\s@endd@\enddocument
\let\bye\relax

\newtoks\tenpoint@ \newtoks\eightpoint@ \newdimen\bigsize@
\newdimen\r@f@nd \newbox\r@f@b@x \newbox\adjb@x \newbox\p@nct@
\newbox\k@yb@x \newcount\rcount
\newbox\b@b@x \newbox\p@p@rb@x \newbox\j@@rb@x \newbox\y@@rb@x
\newbox\v@lb@x \newbox\is@b@x \newbox\p@g@b@x \newif\ifp@g@ \newif\ifp@g@s
\newbox\inb@@kb@x \newbox\b@@kb@x \newbox\p@blb@x \newbox\p@bl@db@x
\newbox\ed@b@x \newif\ifed@ \newif\ifed@s
\newbox\p@p@nf@b@x \newbox\inf@b@x \newbox\b@@nf@b@x

\@ft@\ifx\csname amsppt.sty\endcsname\relax
\headline={\hfil}
\footline={\ifnum\pageno=1\hfil\else\hfil\foliorm\folio\hfil\fi}
\parindent1pc

\font@\tensmc=cmcsc10
\font@\sevenex=cmex7
\font@\sevenit=cmti7
\font@\eightrm=cmr8
\font@\sixrm=cmr6
\font@\eighti=cmmi8 \skewchar\eighti='177
\font@\sixi=cmmi6 \skewchar\sixi='177
\font@\eightsy=cmsy8 \skewchar\eightsy='60
\font@\sixsy=cmsy6 \skewchar\sixsy='60
\font@\eightex=cmex8
\font@\eightbf=cmbx8
\font@\sixbf=cmbx6
\font@\eightit=cmti8
\font@\eightsl=cmsl8
\font@\eightsmc=cmcsc8
\font@\eighttt=cmtt8
\font@\ninerm=cmr9
\font@\ninei=cmmi9 \skewchar\ninei='177
\font@\ninesy=cmsy9 \skewchar\ninesy='60
\font@\nineex=cmex9
\font@\ninebf=cmbx9
\font@\nineit=cmti9
\font@\ninesl=cmsl9
\font@\ninesmc=cmcsc9
\font@\ninemsa=msam9
\font@\ninemsb=msbm9
\font@\nineeufm=eufm9

\UseAMSsymbols \loadeufm
\def\footnoterule{\kern-3\p@\hrule width5pc\kern 2.6\p@}
\def\m@k@foot#1{\insert\footins
 {\interlinepenalty\interfootnotelinepenalty
 \eightpoint\splittopskip\ht\strutbox\splitmaxdepth\dp\strutbox
 \floatingpenalty\@MM\leftskip\z@\rightskip\z@
 \spaceskip\z@\xspaceskip\z@
 \leavevmode\footstrut\ignore#1\unskip\lower\dp\strutbox
 \vbox to\dp\strutbox{}}}
\def\ftext#1{\m@k@foot{\vsk-.8>\nt #1}}
\def\pr@cl@@m#1{\p@sk@p{-100}\medskipamount\b@gr@\nt\ignore
 \bf #1\unskip.\enspace\sl\ignore}
\outer\def\proclaim{\pr@cl@@m} \s@v@d@f\proclaim \let\proclaim\relax
\def\endproclaim{\endgroup\p@sk@p{55}\medskipamount}
\def\demo#1{\sk@@p\medskipamount\nt{\ignore\it #1\unskip.}\enspace
 \ignore}
\def\enddemo{\sk@@p\medskipamount}

\def\cite#1{{\rm[#1]}}
 \def\Refs#1#2{\relax}

\def\big@#1#2{{\hbox{$\left#2\vcenter to#1\bigsize@{}%
  \right.\nulldelimiterspace\z@\m@th$}}}
\def\big{\big@\@ne}
\def\Big{\big@{1.5}}
\def\bigg{\big@\tw@}
\def\Bigg{\big@{2.5}}
\normallineskiplimit\p@

\def\tenpoint{\normalbaselineskip12\p@
 \abovedisplayskip12\p@ plus3\p@ minus9\p@
 \belowdisplayskip\abovedisplayskip
 \abovedisplayshortskip\z@ plus3\p@
 \belowdisplayshortskip7\p@ plus3\p@ minus4\p@
 \textonlyfont@\rm\tenrm \textonlyfont@\it\tenit
 \textonlyfont@\sl\tensl \textonlyfont@\bf\tenbf
 \textonlyfont@\smc\tensmc \textonlyfont@\tt\tentt
 \ifsyntax@ \def\big##1{{\hbox{$\left##1\right.$}}}%
  \let\Big\big \let\bigg\big \let\Bigg\big
 \else
   \textfont\z@\tenrm  \scriptfont\z@\sevenrm
       \scriptscriptfont\z@\fiverm
   \textfont\@ne\teni  \scriptfont\@ne\seveni
       \scriptscriptfont\@ne\fivei
   \textfont\tw@\tensy \scriptfont\tw@\sevensy
       \scriptscriptfont\tw@\fivesy
   \textfont\thr@@\tenex \scriptfont\thr@@\sevenex
        \scriptscriptfont\thr@@\sevenex
   \textfont\itfam\tenit \scriptfont\itfam\sevenit
        \scriptscriptfont\itfam\sevenit
   \textfont\bffam\tenbf \scriptfont\bffam\sevenbf
        \scriptscriptfont\bffam\fivebf
   \setbox\strutbox\hbox{\vrule height8.5\p@ depth3.5\p@ width\z@}%
   \setbox\strutbox@\hbox{\lower.5\normallineskiplimit\vbox{%
        \kern-\normallineskiplimit\copy\strutbox}}%
   \setbox\z@\vbox{\hbox{$($}\kern\z@}\bigsize@1.2\ht\z@
  \fi
  \normalbaselines\rm\dotsspace@1.5mu\ex@.2326ex\jot3\ex@
  \the\tenpoint@}
\def\eightpoint{\normalbaselineskip10\p@
 \abovedisplayskip10\p@ plus2.4\p@ minus7.2\p@
 \belowdisplayskip\abovedisplayskip
 \abovedisplayshortskip\z@ plus2.4\p@
 \belowdisplayshortskip5.6\p@ plus2.4\p@ minus3.2\p@
 \textonlyfont@\rm\eightrm \textonlyfont@\it\eightit
 \textonlyfont@\sl\eightsl \textonlyfont@\bf\eightbf
 \textonlyfont@\smc\eightsmc \textonlyfont@\tt\eighttt
 \ifsyntax@\def\big##1{{\hbox{$\left##1\right.$}}}%
  \let\Big\big \let\bigg\big \let\Bigg\big
 \else
  \textfont\z@\eightrm \scriptfont\z@\sixrm
       \scriptscriptfont\z@\fiverm
  \textfont\@ne\eighti \scriptfont\@ne\sixi
       \scriptscriptfont\@ne\fivei
  \textfont\tw@\eightsy \scriptfont\tw@\sixsy
       \scriptscriptfont\tw@\fivesy
  \textfont\thr@@\eightex \scriptfont\thr@@\sevenex
   \scriptscriptfont\thr@@\sevenex
  \textfont\itfam\eightit \scriptfont\itfam\sevenit
   \scriptscriptfont\itfam\sevenit
  \textfont\bffam\eightbf \scriptfont\bffam\sixbf
   \scriptscriptfont\bffam\fivebf
 \setbox\strutbox\hbox{\vrule height7\p@ depth3\p@ width\z@}%
 \setbox\strutbox@\hbox{\raise.5\normallineskiplimit\vbox{%
   \kern-\normallineskiplimit\copy\strutbox}}%
 \setbox\z@\vbox{\hbox{$($}\kern\z@}\bigsize@1.2\ht\z@
 \fi
 \normalbaselines\eightrm\dotsspace@1.5mu\ex@.2326ex\jot3\ex@
 \the\eightpoint@}

\def\myR@fs{\m@th\s@ct{}{References}\B@gr@\frenchspacing\rcount\z@\widest{AZ}
 \let\Key\key
 \let\refin\relax}
\def\widest#1{\setbox\z@\hbox{[#1]\enspace}\r@f@nd\wd\z@}
\def\R@fb@x{\global\setbox\r@f@b@x} \def\K@yb@x{\global\setbox\k@yb@x}
\def\ref{\par\b@gr@\rm\R@fb@x\box\voidb@x\K@yb@x\box\voidb@x\b@g@nr@f}
\def\c@nc@t#1{\setbox\z@\lastbox
 \setbox\adjb@x\hbox{\unhbox\adjb@x\unhbox\z@\unskip\unskip\unpenalty#1}}
\def\adjust#1{\relax\ifmmode\penalty-\@M\null\hfil$\clubpenalty\z@
 \widowpenalty\z@\interlinepenalty\z@\offinterlineskip\endgraf
 \setbox\z@\lastbox\unskip\unpenalty\c@nc@t{#1}\nt$\hfil\penalty-\@M
 \else\endgraf\c@nc@t{#1}\nt\fi}
\def\adjustnext#1{\P@nct\hbox{#1}\ignore}
\def\cl@s@{\endgraf\setbox\z@\lastbox\global\setbox\@ne\hbox{\unhbox\adjb@x
 \ifvoid\z@\else\unhbox\z@\unskip\unskip\unpenalty\fi}\egroup
 \ifnum\c@rr@nt=\k@yb@x\global\fi
 \setbox\c@rr@nt\hbox{\unhbox\@ne\box\p@nct@}\P@nct\null}
\def\@p@n#1{\def\c@rr@nt{#1}\setbox\c@rr@nt\vbox\bgroup\hsize\maxdimen\nt}
\def\b@g@nr@f{\bgroup\@p@n\z@}
\def\key{\cl@s@\ifvoid\k@yb@x\@p@n\k@yb@x\else\@p@n\z@\fi}
\def\no{\cl@s@\ifvoid\k@yb@x\gad\rcount\n@mbtrue
 \K@yb@x\hbox{\the\rcount}\fi\@p@n\z@}
\def\by{\cl@s@\@p@n\b@b@x} \def\paper{\cl@s@\@p@n\p@p@rb@x\it\ignore}
\def\jour{\cl@s@\@p@n\j@@rb@x} \def\yr{\cl@s@\@p@n\y@@rb@x}
\def\vol{\cl@s@\@p@n\v@lb@x\bf\ignore} \def\issue{\cl@s@\@p@n\is@b@x}
\def\page{\cl@s@\ifp@g@s\@p@n\z@\else\p@g@true\@p@n\p@g@b@x\fi}
\def\pages{\cl@s@\ifp@g@\@p@n\z@\else\p@g@strue\@p@n\p@g@b@x\fi}
\def\inbook{\cl@s@\@p@n\inb@@kb@x} \def\book{\cl@s@\@p@n\b@@kb@x\it\ignore}
\def\publ{\cl@s@\@p@n\p@blb@x} \def\publaddr{\cl@s@\@p@n\p@bl@db@x}
\def\ed{\cl@s@\ifed@s\@p@n\z@\else\ed@true\@p@n\ed@b@x\fi}
\def\eds{\cl@s@\ifed@\@p@n\z@\else\ed@strue\@p@n\ed@b@x\fi}
\def\info{\cl@s@\@p@n\inf@b@x} \def\paperinfo{\cl@s@\@p@n\p@p@nf@b@x}
\def\bookinfo{\cl@s@\@p@n\b@@nf@b@x} \let\finalinfo\info
\def\P@nct{\global\setbox\p@nct@} \def\nopunct{\P@nct\box\voidb@x}
\def\p@@@t#1#2{\ifvoid\p@nct@\else#1\unhbox\p@nct@#2\fi}
\def\sp@@{\penalty-50 \space\hskip\z@ plus.1em}
\def\c@mm@{\p@@@t,\sp@@} \def\sp@c@{\p@@@t\empty\sp@@} \def\p@@nt{.\kern.3em}
\def\p@tb@x#1#2{\ifvoid#1\else#2\@nb@x#1\fi}
\def\@nb@x#1{\unhbox#1\P@nct\lastbox}
\def\endr@f@{\cl@s@\nopunct
 \R@fb@x\hbox{\unhbox\r@f@b@x \p@tb@x\b@b@x\empty
 \ifvoid\j@@rb@x\ifvoid\inb@@kb@x\ifvoid\p@p@rb@x\ifvoid\b@@kb@x
  \ifvoid\p@p@nf@b@x\ifvoid\b@@nf@b@x
  \p@tb@x\v@lb@x\c@mm@ \ifvoid\y@@rb@x\else\sp@c@(\@nb@x\y@@rb@x)\fi
  \p@tb@x\is@b@x{\c@mm@ no\p@@nt}\p@tb@x\p@g@b@x\c@mm@ \p@tb@x\inf@b@x\c@mm@
  \else\p@tb@x \b@@nf@b@x\c@mm@ \p@tb@x\v@lb@x\c@mm@
  \p@tb@x\is@b@x{\sp@c@ no\p@@nt}%
  \ifvoid\ed@b@x\else\sp@c@(\@nb@x\ed@b@x,\space\ifed@ ed.\else eds.\fi)\fi
  \p@tb@x\p@blb@x\c@mm@ \p@tb@x\p@bl@db@x\c@mm@ \p@tb@x\y@@rb@x\c@mm@
  \p@tb@x\p@g@b@x{\c@mm@\ifp@g@ p\p@@nt\else pp\p@@nt\fi}%
  \p@tb@x\inf@b@x\c@mm@\fi
  \else \p@tb@x\p@p@nf@b@x\c@mm@ \p@tb@x\v@lb@x\c@mm@
  \ifvoid\y@@rb@x\else\sp@c@(\@nb@x\y@@rb@x)\fi
  \p@tb@x\is@b@x{\c@mm@ no\p@@nt}\p@tb@x\p@g@b@x\c@mm@ \p@tb@x\inf@b@x\c@mm@\fi
  \else \p@tb@x\b@@kb@x\c@mm@
  \p@tb@x\b@@nf@b@x\c@mm@ \p@tb@x\p@blb@x\c@mm@
  \p@tb@x\p@bl@db@x\c@mm@ \p@tb@x\y@@rb@x\c@mm@
  \ifvoid\p@g@b@x\else\c@mm@\@nb@x\p@g@b@x p\fi \p@tb@x\inf@b@x\c@mm@ \fi
  \else \c@mm@\@nb@x\p@p@rb@x\ic@\p@tb@x\p@p@nf@b@x\c@mm@
  \p@tb@x\v@lb@x\sp@c@ \ifvoid\y@@rb@x\else\sp@c@(\@nb@x\y@@rb@x)\fi
  \p@tb@x\is@b@x{\c@mm@ no\p@@nt}\p@tb@x\p@g@b@x\c@mm@\p@tb@x\inf@b@x\c@mm@\fi
  \else \p@tb@x\p@p@rb@x\c@mm@\ic@\p@tb@x\p@p@nf@b@x\c@mm@
  \c@mm@\@nb@x\inb@@kb@x \p@tb@x\b@@nf@b@x\c@mm@ \p@tb@x\v@lb@x\sp@c@
  \p@tb@x\is@b@x{\sp@c@ no\p@@nt}%
  \ifvoid\ed@b@x\else\sp@c@(\@nb@x\ed@b@x,\space\ifed@ ed.\else eds.\fi)\fi
  \p@tb@x\p@blb@x\c@mm@ \p@tb@x\p@bl@db@x\c@mm@ \p@tb@x\y@@rb@x\c@mm@
  \p@tb@x\p@g@b@x{\c@mm@\ifp@g@ p\p@@nt\else pp\p@@nt\fi}%
  \p@tb@x\inf@b@x\c@mm@\fi
  \else\p@tb@x\p@p@rb@x\c@mm@\ic@\p@tb@x\p@p@nf@b@x\c@mm@\p@tb@x\j@@rb@x\c@mm@
  \p@tb@x\v@lb@x\sp@c@ \ifvoid\y@@rb@x\else\sp@c@(\@nb@x\y@@rb@x)\fi
  \p@tb@x\is@b@x{\c@mm@ no\p@@nt}\p@tb@x\p@g@b@x\c@mm@ \p@tb@x\inf@b@x\c@mm@
 \fi}}
\def\m@r@f#1#2{\endr@f@\ifvoid\p@nct@\else\R@fb@x\hbox{\unhbox\r@f@b@x
 #1\unhbox\p@nct@\penalty-200\enskip#2}\fi\egroup\b@g@nr@f}
\def\endref{\endr@f@\ifvoid\p@nct@\else\R@fb@x\hbox{\unhbox\r@f@b@x.}\fi
 \parindent\r@f@nd\hang
 \textindent{\ifvoid\k@yb@x\else[\unhbox\k@yb@x]\fi}\unhbox\r@f@b@x
 \endgraf\egroup\endgroup}
\def\moreref{\m@r@f;\empty}
\def\transl{\m@r@f;{\unskip\space
 {\sl English translation\ic@}:\penalty-66 \space}}
\def\endRefs{\endgraf\goodbreak\endgroup}

\hyphenation{acad-e-my acad-e-mies af-ter-thought anom-aly anom-alies
an-ti-deriv-a-tive an-tin-o-my an-tin-o-mies apoth-e-o-ses
apoth-e-o-sis ap-pen-dix ar-che-typ-al as-sign-a-ble as-sist-ant-ship
as-ymp-tot-ic asyn-chro-nous at-trib-uted at-trib-ut-able bank-rupt
bank-rupt-cy bi-dif-fer-en-tial blue-print busier busiest
cat-a-stroph-ic cat-a-stroph-i-cally con-gress cross-hatched data-base
de-fin-i-tive de-riv-a-tive dis-trib-ute dri-ver dri-vers eco-nom-ics
econ-o-mist elit-ist equi-vari-ant ex-quis-ite ex-tra-or-di-nary
flow-chart for-mi-da-ble forth-right friv-o-lous ge-o-des-ic
ge-o-det-ic geo-met-ric griev-ance griev-ous griev-ous-ly
hexa-dec-i-mal ho-lo-no-my ho-mo-thetic ideals idio-syn-crasy
in-fin-ite-ly in-fin-i-tes-i-mal ir-rev-o-ca-ble key-stroke
lam-en-ta-ble light-weight mal-a-prop-ism man-u-script mar-gin-al
meta-bol-ic me-tab-o-lism meta-lan-guage me-trop-o-lis
met-ro-pol-i-tan mi-nut-est mol-e-cule mono-chrome mono-pole
mo-nop-oly mono-spline mo-not-o-nous mul-ti-fac-eted mul-ti-plic-able
non-euclid-ean non-iso-mor-phic non-smooth par-a-digm par-a-bol-ic
pa-rab-o-loid pa-ram-e-trize para-mount pen-ta-gon phe-nom-e-non
post-script pre-am-ble pro-ce-dur-al pro-hib-i-tive pro-hib-i-tive-ly
pseu-do-dif-fer-en-tial pseu-do-fi-nite pseu-do-nym qua-drat-ic
quad-ra-ture qua-si-smooth qua-si-sta-tion-ary qua-si-tri-an-gu-lar
quin-tes-sence quin-tes-sen-tial re-arrange-ment rec-tan-gle
ret-ri-bu-tion retro-fit retro-fit-ted right-eous right-eous-ness
ro-bot ro-bot-ics sched-ul-ing se-mes-ter semi-def-i-nite
semi-ho-mo-thet-ic set-up se-vere-ly side-step sov-er-eign spe-cious
spher-oid spher-oid-al star-tling star-tling-ly sta-tis-tics
sto-chas-tic straight-est strange-ness strat-a-gem strong-hold
sum-ma-ble symp-to-matic syn-chro-nous topo-graph-i-cal tra-vers-a-ble
tra-ver-sal tra-ver-sals treach-ery turn-around un-at-tached
un-err-ing-ly white-space wide-spread wing-spread wretch-ed
wretch-ed-ly Brown-ian Eng-lish Euler-ian Feb-ru-ary Gauss-ian
Grothen-dieck Hamil-ton-ian Her-mit-ian Jan-u-ary Japan-ese Kor-te-weg
Le-gendre Lip-schitz Lip-schitz-ian Mar-kov-ian Noe-ther-ian
No-vem-ber Rie-mann-ian Schwarz-schild Sep-tem-ber}
\else
\def\ftext#1{\footnotetext""{\vsk-.8>\nt #1}}

\def\myR@fs{\Refs\nofrills{}\m@th\tenpoint\sect{References}
 \def\k@yf@##1{\hss[##1]\enspace} \let\keyformat\k@yf@
 \def\widest##1{\setbox\z@\hbox{\tenpoint\k@yf@{##1}}\refindentwd\wd\z@}
 \let\Key\key
 \def\refin{\kern\refindentwd}}
\let\info\finalinfo \r@R@fs\Refs
\def\adjust#1{#1} \let\adjustnext\adjust
\fi

\outer\def\myRefs{\myR@fs} \r@st@re\proclaim \r@endd@\bye
\let\Cite\cite \let\Key\key \let\endpro\endproclaim
\let\d@c@\document \def\document{\d@c@\tenpoint}
\hyphenation{ortho-gon-al}

\def\hcor#1{\advance\hoffset by #1}
\def\vcor#1{\advance\voffset by #1}
\let\bls\baselineskip \let\ignore\ignorespaces
\ifx\ic@\undefined \let\ic@\/\fi
\def\vsk#1>{\vskip#1\bls} \let\adv\advance
\def\vv#1>{\vadjust{\vsk#1>}\ignore} \def\vvv#1>{\vadjust{\vskip#1}\ignore}
\def\vvn#1>{\vadjust{\nobreak\vsk#1>\nobreak}\ignore}
\def\vvvn#1>{\vadjust{\nobreak\vskip#1\nobreak}\ignore}
\def\setnormalbls{\edef\normalbls{\bls\the\bls}}
\def\setmaths{\edef\maths{\mathsurround\the\mathsurround}}
\def\Par{\par\medskip} \def\setparindent{\edef\Parindent{\the\parindent}}
\def\Type{\vsk.5>\bgroup\parindent\z@\raggedright\tt}
\def\endType{\vsk.5>\egroup\nt} 

\let\sss\scriptscriptstyle 
\let\vp\vphantom \let\hp\hphantom \let\nt\noindent
  \let\rline\rightline
\def\nn#1>{\noalign{\vskip#1\p@}} \def\NN#1>{\openup#1\p@}
\def\Cup{\bigcup\limits} 
\let\Lim\lim \def\lim{\Lim\limits}
\let\Sum\sum \def\sum{\Sum\limits}
 
\let\Prod\prod \def\prod{\Prod\limits} \let\Int\int \def\int{\Int\limits}

\def\tsum{\mathop{\tsize\Sum}\limits} 

\def\&{.\kern.1em} \def\>{\!\;} \def\){\!\!\;}
\def\nl{\leavevmode\hfill\break}
\def\~{\leavevmode
 \def\n@xt{\ifx\t@st~\def\n@@xt####1{\raise.16ex\mbox{-}}%
\else\def\n@@xt{\kern.04em\raise.16ex\hbox{--}}\fi\n@@xt}\futurelet\t@st\n@xt}

\let\=\m@th \def\mbox#1{\hbox{\m@th$#1$}}
\def\mtext#1{\text{\m@th$#1$}} \def\^#1{\text{\m@th#1}}
\def\Ll@p#1{\llap{\m@th$#1$}} \def\Rl@p#1{\rlap{\m@th$#1$}}
 \def\Cl@p#1{\llap{\m@th$#1$\hss}}
\def\Llap#1{\mathchoice{\Ll@p{\dsize#1}}{\Ll@p{\tsize#1}}{\Ll@p{\ssize#1}}%
 {\Ll@p{\sssize#1}}}
\def\Clap#1{\mathchoice{\Cl@p{\dsize#1}}{\Cl@p{\tsize#1}}{\Cl@p{\ssize#1}}%
 {\Cl@p{\sssize#1}}}
\def\Rlap#1{\mathchoice{\Rl@p{\dsize#1}}{\Rl@p{\tsize#1}}{\Rl@p{\ssize#1}}%
 {\Rl@p{\sssize#1}}}

\let\alb\allowbreak 
\def\ald{\noalign{\alb}} 

\let\o\circ \let\x\times \let\ox\otimes 
\let\sub\subset 
\let\le\leqslant \let\ge\geqslant
\let\der\partial \let\8\infty
\let\bra\langle \let\ket\rangle
 \let\To\Rightarrow
\let\map\mapsto 
 \def\vert{\ |\ } \def\nin{\not\in}
 \let\Empty\varnothing 
 \def\_#1{_{\rlap{$\ssize#1$}}}

\def\lsym#1{#1\alb\ldots\relax#1\alb}
\def\lc{\lsym,}  \def\lx{\lsym\x} \def\lox{\lsym\ox}
\def\llc{\,,\alb\ {\ldots\ ,}\alb\ }
\def\Re{\mathop{\text{\rm Re}\>}} 
\def\E(#1){\mathop{\text{\rm End}\,}(#1)} 
 
\def\Res{\mathop{\text{\rm Res}\>}\limits}

\def\id{\text{\rm id}}  \def\for{\text{for }\,}
\def\1{^{-1}} \def\q{{q\1}}
\def\vst#1{{\lower1.9\p@\mbox{\bigr|_{\raise.5\p@\mbox{\ssize#1}}}}}
\def\qqq{\qquad\quad} 
\def\slap#1{\hbox to\z@{\hss#1\hss}}
\def\sscr#1{\raise.3ex\mbox{\sss#1}}

\def\pms{\raise.25ex\mbox{\ssize\pm}\>}
\def\mps{\raise.25ex\mbox{\ssize\mp}\>}
 \def\mss{{\sscr-}}

 \let\Gm\Gamma
\let\dl\delta \let\Dl\Delta
 \let\eps\varepsilon \let\epsilon\eps
\let\io\iota
\let\ka\kappa
\let\la\lambda \let\La\Lambda
\let\si\sigma  
 \let\phi\varphi
 
\let\thi\vartheta

\let\ze\zeta

\def\C{\Bbb C}
\def\R{\Bbb R}
\def\Z{\Bbb Z}

\def\Zp{\Z_{\ge 0}} \def\Zpp{\Z_{>0}}

\def\OO{\Bbb O}

\def\S{{\bold S}}

\def\Cc{\Cal C}
\def\O{\Cal O}
\def\F{\Cal F}

\def\Sc{\Cal S}
\def\T{\Cal T}

\def\g{\frak g}
\def\CC{\frak C}
\def\ZZ{\frak Z}

\def\TT{\text{\rm T}}

\def\f{\tilde f}
\def\s{\tilde s}
\def\t{\tilde t}

\def\Cbar{\rlap{\=$\kern1.2\p@\overline{\vp{\C}\kern5\p@}$}\C}
\def\Co{\C_{\sss\o}}
\def\Coo{\smash{\rlap{\=$\C_{\sss\o}$}}\phantom{\C}}
\def\CCo{\CC_{\sss\o}} 
\def\CCp{\CC^{\sscr\o}}

\def\Dlp{\Dl^{\!\)\sscr+}} \def\Dlm{\Dl^{\!\)\sscr-}}
\def\Fo{\F_{\!\sss\ox}}
\def\Qo{Q_{\sss\ox}}

\def\Sl{\S_\ell}

\let\M M
\def\gg{\frak{sl}_2} \def\UU{U_q(\gg)}
\def\U{U(\g)} \def\Uq{U_q(\g)}
\def\Un{\U^{\ox n}} \def\Uqn{\Uq^{\ox n}}

\def\v{v^\ast} \def\V{V^\ast}
\def\Vl{V_{[\ell]}} 

\def\Cu{\C\>[u]}
\def\Cn{\C^n} \def\Cl{\C^{\,\ell}}
\def\Con{\Coo^n} \def\Col{\Coo^{\,\ell}}

\def\Zl{\Cal Z_\ell} \def\Zlo{\Zl^{\sscr\o}}

\def\dd#1{{\dsize{\der\over\der#1}}} \def\ddt{\dd{t_a}}
 
\def\ts{t^{\sss\star}}  \def\tmi{t^\mss}
\def\sing{\text{\rm Sing\>}V} \def\singl{\text{\rm Sing\>}\Vl}
\def\dims{\dim\,\singl}

\def\difl/{differential} \def\dif/{difference}
\def\egv/{eigenvector} \def\eva/{eigenvalue} \def\eq/{equation}
\def\lhs/{the left hand side} \def\rhs/{the right hand side}
\def\Lhs/{the left hand side} \def\Rhs/{the right hand side}
\def\gb/{generated by} \def\wrt/{with respect to} \def\st/{such that}
\def\resp/{respectively} \def\off/{offdiagonal} \def\wt/{weight}
\def\pol/{polynomial} \def\rat/{rational} \def\tri/{trigonometric}
\def\fn/{function} \def\var/{variable} \def\raf/{\rat/ \fn/}
\def\sym/{symmetric} \def\perm/{permutation} \def\fd/{finite-dimensional}
\def\rep/{representation} \def\irr/{irreducible} \def\irp/{\irr/ \rep/}
\def\hom/{homomorphism} \def\aut/{automorphism} \def\mult/{multiplicity}
\def\lex/{lexicographical} \def\conv/{convenient} \def\as/{asymptotic}
\def\ndeg/{nondegenerate} \def\neib/{neighbourhood} \def\asex/{\as/ expansion}
\def\deq/{\dif/ \eq/} \def\hw/{highest \wt/} \def\inrp/{integral \rep/}
\def\cc/{compatibility condition} \def\gv/{generating vector}
\def\cp/{critical point} \def\msd/{method of steepest descend}
\def\wlg/{without loss of generality} \def\Wlg/{Without loss of generality}
\def\onedim/{one-dimensional} \def\qcl/{quasiclassical} \def\eqv/{equivalent}
\def\pd/{pairwise distinct} \def\pdj/{pairwise disjunctive}
\def\strh/{string hypothesis} \def\wsep/{well separated}
\def\loc/{local problem} \def\glob/{global problem}
\def\rsp/{regular singular point}

\def\Rm/{\^{$R$-}matrix} \def\Rms/{\^{$R$-}matrices} \def\YB/{Yang-Baxter \eq/}
\def\Ba/{Bethe ansatz} \def\Bv/{Bethe vector} \def\Bae/{\Ba/ \eq/}
\def\KZv/{Knizh\-nik-Zamo\-lod\-chi\-kov}
\def\KZ/{{\sl KZ\/}} \def\qKZ/{{\sl qKZ\/}} \def\KZo/{\qKZ/ operator}

\def\hwm/{\hw/ \gmod/} \def\hwu/{\hw/ \Umod/}
\def\gmod/{\^{\,$\g$-}module} \def\Umod/{\^{$\Uq$-}module}
\def\ggmod/{\^{\,$\gg$-}module}

\def\phf/{phase \fn/} \def\wtf/{\wt/ \fn/} \def\ncp/{\ndeg/ \cp/}
\def\symg/{\sym/ group} \def\sol/{solution} \def\adm/{admissible}
\def\arf/{\adm/ \raf/} \def\mphf/{modified \phf/}
\def\Vval/{\^{$\V\!\!\;$-}valued} \def\Vlv/{\^{$\Vl$-}valued}
\def\Slorb/{\^{$\Sl$-}orbit}

\def\TFT/{Research Insitute for Theoretical Physics}
\def\HY/{University of Helsinki} \def\AoF/{the Academy of Finland}
\def\oldaddress/{P.O\&Box 9 (Siltavuorenpenger 20\,\,C), SF\~~00014, \HY/,
 Finland}
\def\newaddress/{Department of Mathematics, Kyoto University, Kyoto 606,
 Japan}
\def\email/{tarasov\@ kusm.kyoto-u.ac.jp}
\def\SPb/{St\&Petersburg}
\def\home/{\SPb/ Branch of Steklov Mathematical Institute}
\def\UNC/{Department of Mathematics, University of North Carolina}
\def\ChH/{Chapel Hill}
\def\avaddress/{\ChH/, NC 27599, USA} \def\avemail/{av\@ math.unc.edu}
\def\grant/{NSF Grant DMS\~~9203929}

\def\Fadd/{L\&D\&Faddeev} \def\Fre/{I\&Frenkel}
\def\Kir/{A\&N\&Kirillov} \def\Kor/{V\)\&E\&Korepin}
\def\MN/{M\&Nazarov}
\def\Resh/{N\&Reshetikhin} \def\Reshy/{N\&\)Yu\&Reshetikhin}
\def\Skl/{E\&K\&Sklyanin} \def\Takh/{L\&A\&Takhtajan}
\def\Varch/{A\&\)Varchenko} \def\Varn/{A\&N\&\)Varchenko}
\def\VT/{V\)\&Tarasov} \def\VoT/{V\)\&O\&Tarasov}

\def\CMP/{Commun.\ Math.\ Phys.{}}
\def\JSM/{J.\ Sov.\ Math.{}}
\def\TMP/{Theor.\ Math.\ Phys.{}}
\def\ZNS/{Zap.\ nauch.\ semin. LOMI}

\setnormalbls \setmaths \setparindent
\csname langl.def\endcsname

\document
\rline{Kyoto-Math 95\~~04}
\vsk>
\center
\=
{\Bbf Bases of Bethe Vectors and
\vsk.25>
Difference Equations with Regular Singular Points}
\vsk1.5>
\VT/$^{\,\star}$ \ and \ \Varch/$^{\,\ast}$
\vsk>
{\it $^\star$\newaddress/
\vsk.4>
$^\ast$\UNC/%
\ifMag\\ \else, \fi
\avaddress/}
\vsk1.5>
{\sl March \,1995}
\endcenter
\vsk2>
\Abstract
We prove that \Bv/s generically form a base in a tensor product of \irr/
\hw/ \^{$\gg$-}modules or \^{$\UU$-}modules. We apply this result to \deq/s
with \rsp/s. We show that if such an \eq/ has local \sol/s at each of
its singular point, then generically it has a \pol/ \sol/.
\endAbs
\vsk2>
\def\dgg{{\tenpoint\hp{$^\star$}}}
\ftext{\={\tenpoint\sl $^\star$E-mail\/{\rm:} \email/}\nl
\dgg Japan Society for the Promotion of Science Research Fellow\nl
\dgg On leave of absence from \home/
\vv.25>\nl
{\tenpoint\sl $^\ast$E-mail\/{\rm:} \avemail/}\nl
\dgg Supported by \grant/}

\sect{\kern\Parindent Introduction}
\Sno=0
The \Ba/ is a large collection of methods in the theory of quantum integrable
models to calculate the spectrum and \egv/s for a certain commutative
subalgebra of the algebra of observables for an integrable model. This
commutative subalgebra includes the Hamiltonian of the model. Its elements
are called integrals of motion or conservation laws of the model.
The most part of the recent
development of the \Ba/ methods is due to the quantum inverse scattering
transform invented by the Leningrad-\SPb/ school of mathematical physics.
The bibliography on the \Ba/ is enormous. We refer a reader to reviews
\Cite{BIK},\,\Cite{F},\,\Cite{FT}.
\par
Usually in the framework of the \Ba/ \eva/s of conservation laws are expressed
as \fn/s of some additional parameters which have to obey a system of nonlinear
\eq/s. This system is called the system of \Bae/s. In the algebraic \Ba/
there is also a remarkable vector-valued \fn/ of the same additional variables.
Its values at \sol/s to the system of \Bae/s are common \egv/s of conservation
laws. These common \egv/s are called the \Bv/s. An important problem is to
show that the number of appropiate \sol/s to the system of \Bae/s is equal to
the dimension of the representation space of the algebra of observables and
the corresponding \Bv/s form a base of this space.
\par
In this paper we solve this problem for generic integrable models associated
to a finite tensor product of \irr/ \hw/ \^{$\gg$-}modules or
\^{$\UU$-}modules.
A similar but weaker result was announced in a recent preprint \Cite{LS},
see the remark after Theorem 4.2.
\par
A traditional method to study the system of \Bae/s is based on the so-called
\strh/ which predicts a certain behaviour of \sol/s. The \strh/ was used
to get many important results in physical models.
Partially, these results were verified by other methods.
\par
The \strh/ motivated deep combinatorial results, see
\Cite{Ki1},\,\Cite{Ki2}. In particular, it was shown
that the \strh/ predicts the correct number of appropriate
\sol/s to the system of \Bae/s. This fact is called the combinatorial
completeness of the \Bv/s.
\par
A disadvantage of the \strh/ is that it is, strictly speaking,
false and fails to predict even a qualitative picture of \sol/s to the \Bae/s
\Cite{EKS}.
\par
The \strh/ needs to be understood better. We expect that some
clarification of the \strh/ could come from the analysis of the
quantized \KZv/ (\qKZ/\/) \eq/s. In \Cite{TV} we showed that \Bv/s are the
first terms of \as/ \sol/s to the \qKZ/ \eq/s. One could expect that there is
a way to count solutions to the \qKZ/ \eq/s (not the \Bv/s) which is similar
to the counting of the \Bv/s in the \strh/.
\Par
The system of \Bae/s has intimate relation to a certain class of \deq/s which
can be called \deq/s with \rsp/s. These \deq/s arise in two-dimensional
exactly solvable models in statistical mechanics \Cite{B} as well as in
separation of \var/s in quantum lattice integrable models \Cite{S2}, \Cite{S3}.
We apply our results on \Bv/s to these \deq/s to show that such a \deq/ has
a \pol/ \sol/ if it has local \sol/s at singular points. We also count the
number of \deq/s having \pol/ \sol/s.
\Par
The paper is organized as follows. In the first section we study general case
related to a tensor product of \irr/ \hw/ \^{$\gg$-}modules. In the second
section we consider the special case, where the \Bv/s are singular vectors for
the standard $\gg$ action in a tensor product of \^{$\gg$-}modules.
The third section contains the application to \deq/s with \rsp/s.
Results in the $\UU$ case are given in the fourth section for generic $q$ and
in the fifth section for $q$ being a root of unity. In the last section we
consider multiplicative \deq/s with \rsp/s.
\Par
The authors thank R\&P\&Langlands for sending his paper and E\&K\&Sklyanin
for discussions. The first author greatly appreciates hospitality of
the \TFT/, \HY/, where he stayed when this paper had been started.

\Sect{\=Bases of \Bv/s in $\gg$-modules}
Consider the Lie algebra $\g=\gg$ with generators $e,f,h$:
$$
[h,e]=e\,,\qqq [h,f]=-f\,,\qqq [e,f]=2h\,.
$$
Let $\M=\E(\C^2)$. Introduce $T(u)\in\M[u]\ox\U$ as follows:
$$
T(u)=\pmatrix u+h& f\\ e& u-h\endpmatrix\,.
$$
Let $\io_m$ be the embedding $\U\to\Un$ as the \^{$m$-th} tensor factor. Let
$\dl:\U\to\Un$ be the canonical embedding which coincide with
$\io_1\lsym+\io_n$ on $\g$.
\par
Let $z=(z_1\lc z_n)\in\Cn$.
Let $T_m(u)=\id\ox\io_m\bigl(T(u)\bigr)\in\M[u]\ox\Un$. Set
$$
\TT(u)=T_1(u-z_1)\>\ldots\,T_n(u-z_n)=
\pmatrix A(u)& B(u)\\ C(u)& D(u)\endpmatrix
$$
where $A(u)$, $B(u)$, $C(u)$ and $D(u)$ are suitable elements in $\Cu\ox\Un$.
\par
Let $\ka\in\C$. Set $\T(u)=A(u)+\ka D(u)$. Coefficients of the \pol/ $\T(u)$
generate a remarkable commutative subalgebra in $\Un$. The case $\ka=1$ is of
special interest because $\T(u)$ commute with $\dl\bigl(\U\bigr)$ in this
case.
\Rem
The construction explained above originates from the theory of quantum
integrable models. $\T(u)$ is called there the {\it transfer-matrix\/}. Its
coefficients are conservation laws (integrals of motion) of the correspoding
lattice model. For more detailed explanation of the quantum inverse scattering
transform and the algebraic \Ba/ as well as for the related bibliography we
refer a reader to \Cite{BIK},\,\Cite{F},\,\Cite{FT}.
\enddemo
Let $V_1\lc V_n$ be \irr/ \hwm/s with \hw/s $\La_1\lc\La_n$, \resp/.
Set $V=V_1\lox V_n$. A problem of the theory of quantum integrable models
is to diagonalize the operators $\T(u)$ in the space $V$. The {\it algebraic
\Ba/\/} is a tool to construct \egv/s and \eva/s of $\T(u)$. One considers
a special \^{$V\!\!\;$-valued} \fn/ $w(t_1\lc t_\ell)$ of auxiliary \var/s
$t=(t_1\lc t_\ell)$ and chooses $t$ in such a way that $w(t)$ becomes
an \egv/. The arising conditions on $t$ are called the \Bae/s and the
corresponding \egv/ is called a \Bv/. In this paper we prove that for generic
$z_1\lc z_n$ and $\ka$ the \Bv/s form a base in $V$.
\par
More precisely, let $v_1\lc v_n$ be \gv/s of \gmod/s $V_1\lc V_n$, \resp/.
Let $\ell\in\Zp$. Let $\Vl\sub V$ be a \wt/ subspace: $\Vl=\{\,v\in V\vert
\dl(h)\cdot v=\bigl(\tsum_{m=1}^n\La_m-\ell\bigr)\>v\,\}$.
Let $t=(t_1\lc t_\ell)\in\Cl$. Set
$$
w(t)=B(t_1)\ldots B(t_\ell)\cdot v_1\lox v_n\,.
$$
$w(t)$ is a \Vlv/ \sym/ \pol/ in \var/s $t_1\lc t_\ell$. Another formula for
$w(t)$ see below in Lemma 1.12. %% !!! %%
\par
For given $z_1\lc z_n$ consider a system of
algebraic equations on \var/s $t_1\lc t_\ell$:
\ifMag\vv->
$$
\kern-10em
\rightline{$\dsize\prod_{m=1}^n(t_a-z_m+\La_m)\,
\prod_{\tsize{b=1\atop b\ne a}}^\ell(t_a-t_b-1)=
\ka\prod_{m=1}^n(t_a-z_m-\La_m)\,
\prod_{\tsize{b=1\atop b\ne a}}^\ell(t_a-t_b+1)\,.$}
\kern-10em
\Tag{Bae}
$$
\else
$$
\prod_{m=1}^n(t_a-z_m+\La_m)\,
\prod_{\tsize{b=1\atop b\ne a}}^\ell(t_a-t_b-1)=
\ka\prod_{m=1}^n(t_a-z_m-\La_m)\,
\prod_{\tsize{b=1\atop b\ne a}}^\ell(t_a-t_b+1)\,.
\Tag{Bae}
$$
\fi
$a=1\lc\ell$. This system is called the system of {\it \Bae/s\/}.
A \sol/ $t$ to system \(Bae) is called {\it\off/\/} if $t_1\lc t_\ell$
are \pd/, and {\it diagonal\/} otherwise.
\par
Set\vv->
\ifMag
$$
\align
\tau(u,t)&{}=\prod_{m=1}^n(u-z_m+\La_m)\,\prod_{a=1}^\ell{u-t_a-1\over u-t_a}+
\Tag{tau}
\\
&{}+\ka\prod_{m=1}^n(u-z_m-\La_m)\,\prod_{a=1}^\ell{u-t_a+1\over u-t_a}\,.
\endalign
$$
\else
$$
\tau(u,t)=\prod_{m=1}^n(u-z_m+\La_m)\,\prod_{a=1}^\ell{u-t_a-1\over u-t_a}+
\ka\prod_{m=1}^n(u-z_m-\La_m)\,\prod_{a=1}^\ell{u-t_a+1\over u-t_a}\,.
\Tag{tau}
$$
\fi
\Th{eigen}
\back\Cite{BIK},\,\Cite{F},\,\Cite{FT}
Let $t_1\lc t_\ell$ be an \off/ \sol/ to system \(Bae). Then
\ifMag\else\nl\fi
$\T(u)\cdot w(t)=\tau(u,t)\>w(t)$.
\endpro
\Rem
Let $t\in\Cl$ be such that its coordinates are \pd/. Then $t$ is a solution
to system \(Bae) iff $\tau(u,t)$ is a \pol/ in $u$.
\enddemo
Define a set $\ZZ$ in $\Cl$ by the \eq/
$$
\prod_{a=1}^\ell\,\Bigl(\,\prod_{m=1}^n (t_a-z_m+\La_m)(t_a-z_m-\La_m)\,
\prod_{\tsize{b=1\atop b\ne a}}^\ell(t_a-t_b-1)\Bigr)=0\,.
$$
A \sol/ $t$ to system \(Bae) is called {\it \adm/\/} if $t\nin\ZZ$ and
{\it un\adm/\/} otherwise. For an \adm/ \off/ \sol/ $t$ the vector $w(t)$ is
called the {\it \Bv/\/}. A \sol/ $t$ to system \(Bae) is called
a {\it trivial\/} \sol/ if $w(t)=0$ and {\it nontrivial\/} otherwise.
\par
System \(Bae) is preserved by the natural action of the \symg/ $\Sl$ on \var/s
$t_1\lc t_\ell$. Therefore, $\Sl$ acts on \sol/s to this system.
Let $\CC$ be the set of \Slorb/s of \adm/ \off/ \sol/s.
\par
Say that $z_1\lc z_n$ are {\it \wsep/\/} if all points
$z_m-\La_m+s$, $s\in\Zp$, $s<2\La_m$ for $2\La_m\in\Zp$, and $z_m+\La_m$
$m=1\lc n$, are \pd/.
\Rem
$z_1\lc z_n$ are \wsep/ if and only if the module $V$ enjoys
the next properties:
\iitem{i)} $V$ is \irr/ \wrt/ the subalgebra \gb/ coefficients of \pol/s
$A(u),\alb\,B(u),\alb\,C(u),\alb\,D(u)$.
\iitem{ii)} The commutative subalgebra \gb/ coefficients of $A(u)$ acts in $V$
in a semisimple way.
\par\nt
This follows from results of \Cite{T},\,\Cite{NT}.
The part ``{\sl only if\/}'' also follows from results of the present paper.
Note in addition that the action in $V$
of the subalgebra \gb/ coefficients of $A(u)$ has simple spectrum.
\enddemo
\Th{ka<>1}
Let $\ka$ be generic. Let $z_1\lc z_n$ be \wsep/. Then
\iitem{a)} All \adm/ \off/ \sol/s to system \(Bae) are \ndeg/.
\iitem{b)} All un\adm/ \off/ \sol/s to system \(Bae) are trivial.
\iitem{c)} $\#\CC=\dim\Vl$ and the corresponding \Bv/s form a base in $\Vl$.
\endpro
\Rem
This Theorem was proved in \Cite{TV} for generic $z_1\lc z_n$, $\La_1\lc\La_n$.
The proof of Theorem \[ka<>1] is similar to the corresponding proof in
\Cite{TV}. To make this paper self-contained we reproduce some lemmas from
\Cite{TV}.
\enddemo
To prove Theorem \[ka<>1] we use the following strategy. Consider the limit
of system \(Bae) as $\ka\to 0$:
$$
\prod_{m=1}^n\bigl(t_a-z_m+\La_m\bigr)\,
\prod_{\tsize{b=1\atop b\ne a}}^\ell(t_a-t_b-1)=0\,,\qqq a=1\lc\ell\,.
\Tag{Bae0}
$$
All \sol/s to this system are isolated. We consider their deformations
for $\ka\ne 0$. We show that for generic $z_1\lc z_n$ and $\ka$, \off/ \sol/s
to system \(Bae) are \ndeg/ and they are deformations of \off/ \sol/s
to system \(Bae0). The coordinates of an \off/ \sol/ to system \(Bae0) form
a union of arithmetic progressions starting at points $z_1\lc z_n$.
The corresponding \sol/ to system \(Bae) is \adm/ iff each of the arithmetic
progression is not too long.
This proves the first part of claim c). To prove the second part of claim c)
we compute explicitly limits of \Bv/s as $\ka\to 0$ and show that they form
a base in $\Vl$. The claim b) follows from Theorem 1.13. %% !!! %%
\Lm{admi}
Let $2\La_m\nin\Zp$, $m=1\lc n$. For $\ka\ne 0$ and \wsep/ $z_1\lc z_n$ all
\sol/s to system \(Bae) are \adm/.
\endpro
\Pf.
The lemma easily follows from direct analysis of system \(Bae).
As an example consider the case $\ell=2$. Take a \sol/ $t\in\ZZ$.
Suppose $t_1=z_m-\La_m$. Then from the first \eq/, $t_2=t_1+1$ and the second
\eq/ cannot be satisfied. If $t_1=t_2+1$, then from the first \eq/,
$t_1=z_m+\La_m$ for some $m$, and the second \eq/ cannot be satisfied.
Similarly, we can start from $t_1=z_m+\La_m$ or $t_1=t_2-1$. All the other
cases can be obtained by the action of the \symg/ $\S_2$.
\epf
\Lm{isol}
Let $\La=\sum_{m=1}^n\La_m$. Assume that at least one of the following
conditions is fulfilled:
\ifMag\ \else\nl\fi
{\rm a) $\ka\ne1$; \ b) $2\La\nin\Zp$; \ c) $2\La<\ell-1$; \ d) $\La>\ell-1$.}
Then all \off/ \sol/s to system \(Bae) are isolated.
\endpro
\Pf.
Assume that there is a nonisolated \off/ \sol/ to system \(Bae). This means
that we have a curve $t(s)$, $s\in\R$, such that $t(s)$ is an \off/ \sol/ to
system \(Bae) for any $s$. Moreover, we can assume that as $s\to+\8$, $t(s)$
tends to infinity in the following way: $t_a(s)\to\8$ if $a>f$ and $t_a(s)$
has a finite limit if $a\le f$. Set
\ifMag
$$
\align
\tau_0(u)={} &\prod_{m=1}^n (u-z_m+\La_m)\,\prod_{a=1}^f
{u-t_a(+\8)-1\over u-t_a(+\8)}+{}
\Tag{tau0}
\\
&\;{}+\ka\prod_{m=1}^n(u-z_m-\La_m)\,
\prod_{a=1}^f{u-t_a(+\8)+1\over u-t_a(+\8)}\,.
\endalign
$$
\else
$$
\kern-10em
\rightline{$\dsize\tau_0(u)=\prod_{m=1}^n(u-z_m+\La_m)\,\prod_{a=1}^f
{u-t_a(+\8)-1\over u-t_a(+\8)}+ \ka\prod_{m=1}^n(u-z_m-\La_m)\,
\prod_{a=1}^f{u-t_a(+\8)+1\over u-t_a(+\8)}\,.$}
\kern-8em
\Tag{tau0}
$$
\fi
It is clear that for any $u$ \st/ $u\ne t_a(+\8)$, $a=1\lc f$,
$$
\tau(u,t(s))\to\tau_0(u)
\Tag{tauto}
$$
as $\ka\to\ka_0$. Since $t(s)$ is an \off/ \sol/ to system \(Bae),
$\tau(u,t(s))$ is a \pol/ in $u$ for any $s$ and
$\tau_0(u)$ is a \pol/ in $u$ as well. Hence, relation \(tauto) is valid for
any $u$ and coefficients of $\tau(u,t(s))$ tend to coefficients of
$\tau_0(u)$ as $s\to+\8$. From \(tau) and \(tau0) we have that
$$
\tau(u,t(s))-\tau_0(u)=(1-\ka)\>(f-\ell)\>u^{n-1}+\ldots\,.
$$
This means that $\ka=1$. Similarly, for $\ka=1$ we have
$$
\tau(u,t(s))-\tau_0(u)=
\bigl(f(f-2\La-1)-\ell(\ell-2\La-1)\bigr)\>u^{n-2}+\ldots\,.
$$
Since $0\le f<\ell$, this means that $f=2\La+1-\ell$, $2\La\in\Zp$ and
$\ell-1\le 2\La\le 2\ell-2$. The lemma is proved.
\epf
\Lm{88off}
Consider \sol/s to system \(Bae) as (multivalued) \fn/s of $\ka$.
 For any \off/ \sol/ to system \(Bae) every its branch remains finite
for any $\ka\ne 1$.
\endpro
\Pf.
Let $t(\ka)$ be an \off/ \sol/ to system \(Bae). Suppose, $t(\ka)$ tends to
infinity, as $\ka\to\ka_0\ne1$. We can assume that $t_a(\ka)\to\8$ if $a>f$
and $t_a(\ka)$ has a finite limit if $a\le f$. Set
\ifMag
$$
\align
\tau_0(u)={} &\prod_{m=1}^n (u-z_m+\La_m)\,\prod_{a=1}^f
{u-t_a(\ka_0)-1\over u-t_a(\ka_0)}+{}
\Tag{tau00}
\\
&\;{}+\ka_0\prod_{m=1}^n(u-z_m-\La_m)\,
\prod_{a=1}^f{u-t_a(\ka_0)+1\over u-t_a(\ka_0)}\,.
\endalign
$$
\else
$$
\kern-10em
\rightline{$\dsize\tau_0(u)=\prod_{m=1}^n(u-z_m+\La_m)\,\prod_{a=1}^f
{u-t_a(\ka_0)-1\over u-t_a(\ka_0)}+ \ka_0\prod_{m=1}^n(u-z_m-\La_m)\,
\prod_{a=1}^f{u-t_a(\ka_0)+1\over u-t_a(\ka_0)}\,.$}
\kern-8em
\Tag{tau00}
$$
\fi
Since $t(\ka)$ is an \off/ \sol/ to system \(Bae), $\tau(u,t(\ka))$ is a \pol/
in $u$ for any $\ka\ne\ka_0$ and $\tau_0(u)$ is a \pol/ in $u$ as well.
Similar to the proof of Lemma \[isol] we have that coefficients of
$\tau(u,t(\ka))$ tend to coefficients of $\tau_0(u)$ as $\ka\to\ka_0$.
On the other hand, it follows from \(tau) and \(tau00) that
$$
\tau(u,t(\ka))-\tau_0(u)=(\ka-\ka_0)\>
\bigl(u^n+\bigl(\ell-\tsum_{m=1}^n(z_m+\La_m)\bigr)\>u^{n-1}\bigr)+
(1-\ka_0)\>(f-\ell)\>u^{n-1}+\ldots
$$
which implies $f=\ell$. The lemma is proved.
\epf
Set
\vv.2>
$\OO=\{\,\eta\in\Z^{n\ell}_{\ge 0}\vert
\sum_{m=1}^n\sum_{a=1}^\ell\eta_{ma}=\ell$ and $\eta_{ma}=0\To\eta_{m,a+1}=0$
$m=1\lc n$, $a=1\lc\ell\,\}$. Let $\O(\eta)$ be an \Slorb/ of \sol/s to system
\(Bae0) fixed by conditions
$$
\#\{b\ |\ t_b=z_m-\La_m+a-1\}=\eta_{ma}\,,\qquad
m=1\lc n\,,\quad a=1\lc \ell\,.
$$
\Lm{off}
\iitem{a)} All \sol/s to system \(Bae0) are isolated.
\iitem{b)} Any \Slorb/ of \sol/ to system \(Bae0) has the form $\O(\eta)$ for
a suitable $\eta\in\OO$.
\iitem{c)} For generic $z_1\lc z_n$ \Slorb/s $\O(\eta)$ of \sol/s to system
\(Bae0) are \pd/.
\iitem{d)} For generic $z_1\lc z_n$ all \off/ \sol/s to system \(Bae0)
are \ndeg/.
\endpro
\Pf.
Claim a) follows from the next lemma.
\Lm{isol1}
Let $Q_1\lc Q_l$ be \pol/s in one \var/, $\deg Q_a>0$, $a=1\lc l$.
Then all \sol/s to the system
$$
Q_a(t_a)\prod_{\tsize{b=1\atop b\ne a}}^l(t_a-t_b-1)=0\,,\qqq a=1\lc l\,,
$$
for \var/s $t_1\lc t_l$ are isolated.
\endpro
\Pf.
We prove the lemma by induction on $l$. The case $l=1$ is clear.
\par
Let $t_1\lc t_l$ be a \sol/ to the system in question. \Wlg/ we can assume that
$\Re t_1\lsym=\Re t_f<\Re t_{f+1}\lsym\le\Re t_l$. Therefore, we have that
$Q_a(t_a)=0$, $a=1\lc f$. Remaining \var/s $t_{f+1}\lc t_l$ satisfies
a new system
$$
\gather
\tilde Q_a(t_a)\prod_{\tsize{b=f+1\atop b\ne a}}^l(t_a-t_b-1)=0\,,
\qqq a=f+1\lc l\,,
\Tag{new1}
\\
\Text{where}
\nn-15>
\tilde Q_a(t_a)=Q_a(t_a)\prod_{b=1}^f(t_a-t_b-1)\,.
\endgather
$$
By the induction assumption all \sol/ to system \(new1) are isolated. Since
possible values of $(t_1\lc t_f)$ form a discrete set, the lemma is proved.
\epf
It is clear that the proof of Lemma \[isol1] implies claim b) as
well. Claim c) is evident. To prove claim d) we drop all factors which remain
nonzero on the \sol/ in question and the apply the next lemma.
\Lm{mult}
Let $Q_1\lc Q_l$ be homogeneous \pol/s in variables $x_1\lc x_l$.
Assume, that $x_a=0$, $a=1\lc l$, is an isolated solution to a system
$$
Q_a(x)=0\,,\qquad a=1\lc l\,.
$$
Then the \mult/ of this solution is equal to $\prod_{a=1}^l\deg Q_a$.
\endpro
\nt
Lemma \[isol1] is proved.
\epf
\Rem
Lemma \[mult] corresponds to Lemma 4.9 in \Cite{TV}. The form in which
Lemma 4.9 is formulated in \Cite{TV} is wrong. References to Lemma 4.9
in \Cite{TV} must be replaced by references to the above Lemma \[mult].
All results in \Cite{TV} remain correct.
\enddemo
\Lm{dd}
Let $z_1\lc z_n$ be generic.
Let $t(\ka)$ be a \sol/ to system \(Bae), which is a deformation of
a diagonal \sol/ $t(0)$ to system \(Bae0). Then $t(\ka)$ is a diagonal \sol/.
\endpro
\Pf.
Let $a_0=\min\,\{\,a\ |\ \eta_{ma}>1$ for some $m\,\}$. Let $m_0$ be \st/
$\eta_{m_0a_0}>1$. Let, for example, $t_1=t_\ell=z_{m_0}+a_0-1$. Let
$\tmi=(t_1\lc t_{\ell-1})\in\C^{\,\ell-1}$.
Consider a new system
\ifMag
$$
\gather
\kern-10em
\rightline{$\dsize
\prod_{m=1}^n(t_1-z_m+\La_m)\,\prod_{b=2}^{\ell-1}(t_1-t_b-1) =
-\ka\prod_{m=1}^n(t_1-z_m-\La_m)\,\prod_{b=2}^{\ell-1}(t_1-t_b+1)\,,$}
\kern-10em
\Tag{new0}
\\
\nn6>
\aligned
(t_a &{}-t_1-1)\prod_{m=1}^n(t_a-z_m+\La_m)\,
\prod_{\tsize{b=1\atop b\ne a}}^{\ell-1}(t_a-t_b-1) ={}
\\
\nn-4>
&{}=\ka\,(t_a-t_1+1)\prod_{m=1}^n(t_a-z_m-\La_m)\,
\prod_{\tsize{b=1\atop b\ne a}}^{\ell-1}(t_a-t_b+1)\,,
\endaligned
\endgather
$$
\else
$$
\align
\prod_{m=1}^n(t_1-z_m+\La_m)\,\prod_{b=2}^{\ell-1}(t_1-t_b-1) &=
-\ka\prod_{m=1}^n(t_1-z_m-\La_m)\,\prod_{b=2}^{\ell-1}(t_1-t_b+1)\,,
\Tag{new0}
\\
\nn3>
(t_a-t_1-1)\prod_{m=1}^n(t_a-z_m+\La_m)\,
\prod_{\tsize{b=1\atop b\ne a}}^{\ell-1}(t_a-t_b-1) &=
\ka\,(t_a-t_1+1)\prod_{m=1}^n(t_a-z_m-\La_m)\,
\prod_{\tsize{b=1\atop b\ne a}}^{\ell-1}(t_a-t_b+1)\,,
\endalign
$$
\fi
$a=2\lc\ell-1$. It is obtained from system \(Bae) by the substitution
$t_\ell=t_1$. Incidently, the first and the last \eq/s of \(Bae)
coincide after this substitution. Any \sol/ to system \(new0) gives rise to
a diagonal \sol/ to system \(Bae) by setting $t_\ell=t_1$.
\par
By Lemma \[isol1] $\tmi(0)$ is an isolated \sol/ to system \(new0) at $\ka=0$.
It follows from Lemma \[mult] that the \sol/ $t(0)$ to system \(Bae0) has
the same \mult/ as the \sol/ $\tmi(0)$ to system \(new0) at $\ka=0$. Namely,
the prescribed choice of coinciding coordinates of the \sol/ $t(0)$ to system
\(Bae0) guarantees that a homogeneous \pol/ corresponding to the dropped \eq/
with $a=\ell$ is of degree $1$. This choice also guarantees that degrees of
homogeneous \pol/s corresponding to the remaining \eq/s with $a=1\lc\ell-1$
are the same for both systems \(Bae0) and \(new0) at $\ka=0$.
\par
The coincidence of multiplicities means that any
deformation of the diagonal \sol/ $t(0)$ can be obtained from
some deformation of the \sol/ $\tmi(0)$ and, therefore, is diagonal.
\epf
\Pf of Theorem \[ka<>1].
Henceforward, until the end of the section we assume that $\ka$ is generic,
unless the contrary is indicated explicitly. Let $z_1\lc z_n$ be generic,
until the contrary is indicated explicitly.
\par
Set $\Zl=\{\,\nu\in\Z^n_{\ge 0}\vert\sum_{m=1}^n\nu_m=\ell\,\}$.
\vv.1>
Set $\Zlo=\{\,\nu\in\Z^n_{\ge 0}\vert\sum_{m=1}^n\nu_m\alb=\ell$ and
$\nu_m\le 2\La_m$, if $2\La_m\in\Zp$, $m=1\lc n\,\}$. For
$\nu\in\Zl$ let $\ts(\nu)$ be the following \off/ \sol/ to system \(Bae0)
$$
\ts_a(\nu)=z_m-\La_m+a-\ell_{m-1}-1\,,\qqq\ell_{m-1}<a\le\ell_m\,,
\Tag{tnu}
$$
where $\ell_m=\sum_{k=1}^m\nu_k$, \,$\ell_0=0$, $\ell_n=\ell$.
It is related to the previous description of \sol/s as $\eta_{ka}=1$ for
$a\le\nu_k$ and $\eta_{ka}=0$ for $a>\nu_k$.
\par
Any \Slorb/ of \off/ solutions to system \(Bae0) can be obtained as the
orbit of a solution $\ts(\nu)$ for a suitable $\nu$. Moreover, for generic
$z_1\lc z_n$ \sol/s $\ts(\nu)$ for different $\nu$ belong to different
\Slorb/s.
Let $t(\nu,\ka)$ be a \sol/ to system \(Bae) which is the deformation of
$\ts(\nu)$.
\Lm{Zlo}
Let $\nu\in\Zlo$. Then $t(\nu,\ka)$ is an \adm/ \sol/.
\endpro
\Pf.
By direct analysis of system \(Bae) we obtain that
$$
t_a(\nu,\ka)=\ts(\nu)+\ka^{\ell_i-a+1}t'_a(\nu,\ka)\,,\qqq
t'_a(\nu,0)\ne0\,,
$$
for $\ell_{i-1}<a\le\ell_i$. This means that $t(\nu,\ka)$ is an \adm/ \sol/
for small $\ka$, and hence, for generic $\ka$.
\epf
\Lm{Zloo}
Let $\nu\in\Zl\setminus\Zlo$. Then $t(\nu,\ka)$ is an un\adm/ \sol/.
\endpro
\Pf.
Let, for example, $\nu_1>2\La_1$. Set $d=2\La_1+1$.
Let $\tmi=(t_{d+1}\lc t_\ell)\in\C^{\,\ell-d}$.
Consider a system for $t_{d+1}\lc t_\ell$
\ifMag
$$
\kern-10em
\rightline{$\dsize\prod_{m=1}^n(t_a-z_m+\La_m)\,
\prod_{\tsize{b=1\atop b\ne a}}^\ell(t_a-t_b-1)=
\ka\prod_{m=1}^n(t_a-z_m-\La_m)\,
\prod_{\tsize{b=1\atop b\ne a}}^\ell(t_a-t_b+1)\,,$}
\kern-10em
\Tag{new}
$$
\else
$$
\prod_{m=1}^n\bigl(t_a-z_m+\La_m\bigr)\,
\prod_{\tsize{b=1\atop b\ne a}}^\ell(t_a-t_b-1)=
\ka\prod_{m=1}^n\bigl(t_a-z_m-\La_m\bigr)\,
\prod_{\tsize{b=1\atop b\ne a}}^\ell(t_a-t_b+1)\,,
\Tag{new}
$$
\fi
$a=d+1\lc\ell$ where $t_c=\ts_c(\nu)$ for $c=1\lc d$. This system is obtained
from system \(Bae) by the substitution $t_c=\ts_c(\nu)$, $c=1\lc d$.
Incidently, the first $d$ equations of \(Bae) become identities after this
substitution. Any \sol/ to system \(new) gives rise to an un\adm/ \sol/ to
system \(Bae) by setting $t_c=\ts_c(\nu)$, $c=1\lc d$. The \mult/ of \sol/
$\ts(\nu)$ to system \(Bae0) is equal to $1$ since $\ts(\nu)$ is \off/. This
means that the \mult/ of the \sol/ $\tmi(\nu,0)$ to system \(new) at $\ka=0$
is also equal to $1$. Moreover, the \sol/ $\t(\nu,\ka)$ to system \(Bae) can
be obtained from the \sol/ $\tmi(\nu,\ka)$ to system \(new) and, therefore,
is un\adm/.
\epf
 Finally, there are precisely $\#\Zlo=\dim\Vl$ \Slorb/s of \adm/ \off/ \sol/s.
All of them are the orbits of \sol/s $t(\nu,\ka)$ for a suitable $\nu\in\Zlo$
and, hence, are \ndeg/.
\Par
Introduce the canonical monomial base in $V$:
$\{\,F^\nu = f^{\nu_1}v_1\lox f^{\nu_n}v_n\,\}$. It is clear that
$\{\,F^\nu\vert\nu\in\Zlo\,\}$ is a base in $\Vl$.
\Lm{decomp}
\back\Cite{Ko}
The following decomposition holds
$$
w(t)=\sum_{\text{\it part}}\,F^{\nu(\Gm)}\,\prod_{l=2}^n\,\prod_{m=1}^{l-1}\,
\Bigl(
\prod_{\tsize{a\in\Gm_l\atop\rlap{\=$\ssize b\in\Gm_m$}\hp{a\in\Gm_l}}}
{t_a-t_b-1\over t_a-t_b}
\prod_{a\in\Gm_l}(t_a-z_m+\La_m)\prod_{a\in\Gm_m}(t_a-z_l-\La_l)\Bigr)\,.
$$
Here the sum is taken over all partitions of the set $\{1\lc\ell\}$ into
disjoint subsets $\Gm_1\lc\Gm_n$ and $\nu(\Gm)=(\#\Gm_1\lc\#\Gm_n)$.
\endpro
Define on $\Zl$ a \lex/ order: $\nu<\nu'$ if $\nu_1<\nu'_1$, or
$\nu_1=\nu'_1$, $\nu_2<\nu'_2$, etc. Say that $b\ll a$ if $b\le\ell_m<a$ for
some $m$. Lemma \[decomp] implies that for $t=\ts(\nu)$
\ifMag
$$
\align
w(t)= F^\nu\prod_{a=2}^\ell\,\prod_{b\ll a}{t_a-t_b-1\over t_a-t_b}\,
\prod_{m=1}^n\,\Bigl(\prod_{a>\ell_m}(t_a-z_m+\La_m)
\prod_{a\le\ell_{m-1}}(t_a-z_m-\La_m)\Bigr)+{}&
\\
\nn4>
{}+\sum_{\nu'>\nu}F^{\nu'}\theta_{\nu\nu'}\,& .
\endalign
$$
\else
$$
w(t)= F^\nu\prod_{a=2}^\ell\,\prod_{b\ll a}{t_a-t_b-1\over t_a-t_b}\,
\prod_{m=1}^n\,\Bigl(\prod_{a>\ell_m}(t_a-z_m+\La_m)
\prod_{a\le\ell_{m-1}}(t_a-z_m-\La_m)\Bigr)
+ \sum_{\nu'>\nu}F^{\nu'}\theta_{\nu\nu'}\,.
$$
\fi
where $\theta_{\nu\nu'}$ are suitable coefficients. This means that
$\{\,w\bigl(t(\nu,\ka)\bigr)\vert\nu\in\Zlo\,\}$ is a base in $\Vl$. Hence,
$\{\,w\bigl(t(\nu,\ka)\bigr)\vert\nu\in\Zlo\,\}$ is a base in $\Vl$
for generic $\ka$.
\Par
Let $\v_m$ be a linear \fn/ on $V_m$ \st/ $\bra\v_m,v_m\ket=1$ and
$\bra\v_m,v\ket=0$ for any \wt/ vector $v\in V_m$, $v\ne v_m$.
\Th{ortho}
\back\Cite{Ko},\,\Cite{TV}
Let $\La_1\lc\La_n$, $z_1\lc z_n$ and $\ka$ be generic.
\iitem{a)} For any \off/ \sol/ $t=(t_1\lc t_n)$ to system \(Bae)
\par
\ifMag
$$
\align
\bra\v_1\lox{}& \v_n,C(t_1)\ldots C(t_n)w(t)\ket={}
\\
\nn3>
&\)\llap{${}={}$}\>(-1)^\ell
\prod_{m=1}^n\,\prod_{a=1}^\ell\bigl(t_a-z_m-\La_m)\,
\prod_{a=2}^\ell\,\prod_{b=1}^{a-1}{t_a-t_b+1\over t_a-t_b}\,\x{}
\\
\nn4>
&\>\llap{${}\x{}$}\)\det\left[
\aligned
\ddt\Bigl(\> &\prod_{m=1}^n(t_b-z_m+\La_m)
\,\prod_{\tsize{c=1\atop c\ne b}}^\ell(t_b-t_c-1)-{}
\\
{}- &\prod_{m=1}^n(t_b-z_m-\La_m)
\,\prod_{\tsize{c=1\atop c\ne b}}^\ell(t_b-t_c+1)\Bigr)
\endaligned
\right]_{a,b=1\lc\ell}\,.
\endalign
$$
\else
$$
\align
\bra &\v_1\lox\v_n,C(t_1)\ldots C(t_n)w(t)\ket=
(-1)^\ell\prod_{m=1}^n\,\prod_{a=1}^\ell(t_a-z_m)-\La_m)\,
\prod_{a=2}^\ell\,\prod_{b=1}^{a-1}{t_a-t_b+1\over t_a-t_b}\>\x{}
\\
&{}\x\>\det\Biggl[\ddt\Bigl(\>\prod_{m=1}^n(t_b-z_m+\La_m)\,
\prod_{\tsize{c=1\atop c\ne b}}^\ell(t_b-t_c-1)-
\ka\prod_{m=1}^n(t_b-z_m-\La_m)\,
\prod_{\tsize{c=1\atop c\ne b}}^\ell(t_b-t_c+1)\Bigr)
\Biggr]_{a,b=1\lc\ell}\,.
\endalign
$$
\fi
\iitem{b)} For any \off/ \sol/s $t$ and $\t=(\t_1\lc\t_n)$ to system \(Bae)
which lie in different \Slorb/s
$$
\bra\v_1\lox\v_n,C(\t_1)\ldots C(\t_n)w(t)\ket=0\,.
$$
\iitem{c)} Let $t,\t$ be isolated \sol/s to system \(Bae). Then both claims
{\rm a)} and {\rm b)} remain valid for any $\La_1\lc\La_n$, $z_1\lc z_n$ and
$\ka$.
\endpro
\Rem
In this paper we use a normalization of $w(t)$ which differs from the
normalization in \Cite{TV}.
\enddemo
Let $\t$ be an un\adm/ \off/ \sol/ to system \(Bae). Since
$\{\,w\bigl(t(\nu,\ka)\bigr)\vert\nu\in\Zlo\,\}$ is a base in $\Vl$ for
generic $\ka$ and $z_1\lc z_n$, Theorem \[ortho] implies $w(\t)$=0.
The theorem is proved for generic $z_1\lc z_n$.
\Par
Let $z_1\lc z_n$ be \wsep/, but not necessarily generic.
\Lm*
Let $\ka$ be generic. Then $t(\nu,\ka)$, $\nu\in\Zlo$, are \ndeg/ \adm/ \off/
\sol/s to system \(Bae) and the corresponding \Bv/s form a base in $\Vl$.
\endpro
\nt
The proof is the same as for generic $z_1\lc z_n$.
\Par
Since by Lemma \[isol] all \off/ \sol/s to system \(Bae) are isolated, they
can be deformed to the case of generic $z_1\lc z_n$. Now the theorem follows
from the last lemma and Theorem \[ka<>1] for generic $z_1\lc z_n$.
\epf
\Rem
There is another proof to claim b) of Theorem \[ka<>1]. We give it below.
\Lm{chain}
Let $\ka$ be generic. Let $z_1\lc z_n$ be \wsep/. For any un\adm/ \off/ \sol/
$t_1\lc t_\ell$ to system \(Bae) the set $\{\>t_1\lc t_\ell\>\}$ contains
a set $\{\>z_m-\La_m,\,z_m-\La_m+1,\lc z_m+\La_m\>\}$ for some $m$.
\endpro
\Pf.
Already established claims a) and c) of Theorem \[ka<>1] imply that it
suffices to prove the lemma only for generic $z_1\lc z_n$. In the last case
any \off/ \sol/ to system \(Bae) belongs to the orbit of \sol/ $t(\nu,\ka)$
for a suitable $\nu$. Hence, the lemma follows from the proof to Lemma \[Zloo].
\epf
\Th{sigma}
\back{\Cite{T} {\rm (cf.\ \Cite{CP1} for more detailed proof)}}
Let $\si$ be a \perm/ of $1\lc n$. Set $V^\si=V_{\si(1)}\lox V_{\si(n)}$
and $z^\si=(z_{\si(1)}\lc z_{\si(n)})$. Let $w^\si(t)$ be constructed in the
same manner for $V^\si, z^\si$ as $w(t)$ is constructed for $V,z$. Then
for \wsep/ $z_1\lc z_n$ there is a linear isomorphism $V\to V^\si$ \st/
$w(t)\map w^\si(t)$.
\endpro
By Lemma \[chain] and taking into account Theorem \[sigma] we can assume that
$\t_1=z_1-\La_1$, $\t_2=z_1-\La_1+1\llc \t_{2\La_1+1}=z_1+\La_1$
\wlg/. The structure of the formula in Lemma \[decomp] is \st/ the scalar
factor is zero unless $\t_m\in\Gm_1$, $m=1\lc 2\La+1$. Hence, only terms with
$\#\Gm_1\ge 2\La_1+1$ survive in the sum. But in these terms $F^{\nu(\Gm)}=0$
since $f^{2\La_1+1}v_1=0$. Therefore, $w(\t)=0$.
\epf

\Sect{\=Bases of \Bv/s in $\gg$-modules; bases of singular vectors}
In this section we always assume that $\ka=1$.
Set $\sing=\{\,v\in V\vert\dl(e)\cdot v=0\,\}$ and $\singl=\Vl\cap\sing$.
Let $\CCo$ be the set of \Slorb/s of nontrivial isolated \adm/ \off/ \sol/s
to system \(Bae).
\Lm{sing}
\back\Cite{F},\,\Cite{FT2}
Let $t$ be an \off/ \sol/ to system \(Bae). Then $w(t)\in\singl$.
\endpro
\Th{ka=1}
Let $\ka=1$. Then
\iitem{a)} For generic $z_1\lc z_n$, all nontrivial isolated \adm/ \off/
\sol/s to system \(Bae) are \ndeg/.
\iitem{b)} For any $z_1\lc z_n$, all trivial \adm/ \off/ \sol/s to system
\(Bae) are degenerate.
\iitem{c)} For generic $z_1\lc z_n$, all isolated un\adm/ \off/ \sol/s
to system \(Bae) are trivial.
\iitem{d)} For generic $z_1\lc z_n$, $\#\CCo=\dims$ and the corresponding
\Bv/s form a base in $\singl$.
\endpro
\Pf.
A trivial \adm/ \sol/ to system \(Bae) is either nonisolated, and hence
degenerate, or degenerate by Theorem \[ortho]. This proves claim b).
\Par
Let $s\in\C$, $s\ne 0$. Make a change of \var/s $x=sz\in\Cn$, $u=st\in\Cl$.
In the new \var/s $u_1\lc u_\ell$ system \(Bae) reads as follows:
\ifMag
$$
\kern-10em
\rightline{$\dsize\prod_{m=1}^n(t_a-z_m+s\La_m)\,
\prod_{\tsize{b=1\atop b\ne a}}^\ell(t_a-t_b-s)=
\ka\prod_{m=1}^n(t_a-z_m-s\La_m)\,
\prod_{\tsize{b=1\atop b\ne a}}^\ell(t_a-t_b+s)\,.$}
\kern-10em
\Tag{Baen}
$$
\else
$$
\prod_{m=1}^n(t_a-z_m+s\La_m)\,
\prod_{\tsize{b=1\atop b\ne a}}^\ell(t_a-t_b-s)=
\ka\prod_{m=1}^n(t_a-z_m-s\La_m)\,
\prod_{\tsize{b=1\atop b\ne a}}^\ell(t_a-t_b+s)\,.
\Tag{Baen}
$$
\fi
$a=1\lc\ell$. As $s\to 0$ system \(Baen) turns into
$$
\prod_{m=1}^n(u_a-x_m)\,\prod_{\tsize{b=1\atop b\ne a}}^\ell(u_a-u_b)\,
\Bigl(\,\sum_{m=1}^n{\La_m\over u_a-x_m}\,-
\sum_{\tsize{b=1\atop b\ne a}}^\ell{1\over u_a-u_b}\,\Bigr)=0\,,
\Tag{Baen0}
$$
$a=1\lc\ell$. Both systems \(Baen) and \(Baen0) are preserved by the
natural action of the \symg/ $\Sl$ on variables $u_1\lc u_\ell$.
\Lm{RV}
\back\Cite{RV}
Let $x_1\lc x_n$ be generic. Then there are $\dims$ \Slorb/s of \ndeg/ \off/
\sol/s to system \(Baen0).
\endpro
Let $\CCp$ be the set of \Slorb/s of \ndeg/ \adm/ \off/ \sol/s to system
\(Bae). Since system \(Bae) is a deformation of system \(Baen0), then
Lemma \[RV] means that $\#\CCp\ge\dims$ for generic $z_1\lc z_n$. On the other
hand Lemma \[sing] and Theorem \[ortho] show $\#\CCp\le\dims$. Then
$\#\CCp=\dims$. Moreover, the corresponding \Bv/s form a base in $\singl$.
\par
Let $\t$ be a solution to system \(Bae) \st/ $\Sl(\t)\in\CCo\setminus\CCp$.
Since $\{\,w(t)\vert\Sl(t)\in\CCp\,\}$ is a base in $\singl$ Theorem \[ortho]
implies that $w(\t)=0$, and hence, $\CCo\sub\CCp$. The opposite inclusion
$\CCp\sub\CCo$ clearly follows from claim b), and we obtain that $\CCo=\CCp$.
\par
If $\t$ is an isolated un\adm/ \off/ \sol/ to system \(Bae) then
Theorem \[ortho] again implies that $w(\t)=0$ since
$\{\,w(t)\vert\Sl(t)\in\CCp\,\}$ is a base in $\singl$. The theorem is proved.
\epf
\Th{details}
Let $\ka=1$. Let $z_1\lc z_n$ be generic.
\iitem{a)} Let $2\La_m\nin\Zp$ for some $m$. Then all degenerate \adm/ \off/
\sol/s to system \(Bae) are nonisolated.
\iitem{b)} Let $2\La_m\in\Zp$, $m=1\lc n$. Then all degenerate or un\adm/
\off/ \sol/s to system \(Bae) are trivial.
\endpro
\Pf.
Assume, for example that $2\La_1\nin\Zp$. All other cases can be considered
similarly. Let $t\nin\ZZ$. By Lemma \[decomp] we have
$$
w(t)=f^\ell v_1\ox v_2\lox v_n\prod_{a=1}^\ell(t_a-z_1-\La_1)\,
+\,\sum_{k=0}^{\ell-1}f^kv_1\ox w_k(t)\,\ne 0\,.
$$
Here $w_k(t)$ are suitable \^{$V_2\lox V_n$-}valued \pol/s. Since by
Theorem \[ka=1] for generic $z_1\lc z_n$ any isolated degenerate \adm/ \off/
\sol/ to system \(Bae) is trivial, claim a) is proved.
\par
By Theorem \[ka=1] any degenerate or un\adm/ \off/ \sol/ to system \(Bae) is
either trivial or nonisolated. By Lemma \[isol] existence of nonisolated \off/
\sol/s implies that $\ell\ge\La+1$. But $\singl=0$ for $\ell>\La$ in the case
in question. Lemma \[sing] completes the proof.
\epf

\Sect{Difference \eq/s with \rsp/s. Additive case}
Consider a second order \difl/ \eq/ with \pol/ coefficients
$$
L(u)\xi''(u)+M(u)\xi'(u)+N(u)\xi(u)=0\,,
\Tag{Gaud}
$$
$\deg L=n$, $\deg M=n-1$.
Assume that this \eq/ has $n$ \pd/ \rsp/s
$z_1\lc z_n$ with exponents $0,\la_m+1$, $\la_m\in\Zpp$, at a point $z_m$,
$m=1\lc n$. This means that $L(u)=\prod_{m=1}^n(u-z_m)$ and
$M(u)/L(u)=-\sum_{m=1}^n\la_m/(u-z_m)$.
\Par\nt
{\it Global problem.}\enspace
To determine a \pol/ $N(u)$ \st/ \eq/ \(Gaud) has a nonzero \pol/ \sol/.
\par\nt
{\it Local problem.}\enspace
To determine a \pol/ $N(u)$, $\deg N\le n-2$, \st/ all \sol/s to \eq/ \(Gaud)
are entire \fn/s.
\Par
These problems arise in separation of \var/s in the Gaudin model \Cite{S1},
see also \Cite{Sz, Sect.\ 6.8}, \Cite{RV}.
The next theorem, first observed by Sklyanin \Cite{S1},
easily follows from analytic theory of \difl/ \eq/s.
\Th*
Let $N(u)$ be a \sol/ to the \loc/. Then $N(u)$ is a \sol/ to the \glob/.
\endpro
\demo{\it Problem}
For fixed $\la_1\lc \la_n$ and $z_1\lc z_n$ to determine the number of \sol/s
to the formulated problems, see \Cite{S1}, \Cite{Sz, Sect.\ 6.8}, \Cite{RV}.
\enddemo
In this section we consider similar problems for \deq/s arising in separation
of \var/s in quantum lattice integrable models \Cite{S2}, \Cite{S3}.
Similar \deq/s were introduced by Baxter in his famous studies of
two-dimensional exactly solvable models in statistical mechanics \Cite{B}.
\Par
We continue to use notations of previous sections. All over this section
we assume that $2\La_m\in\Zp$, $m=1\lc n$. We also assume that $z_1\lc z_n$
are \wsep/, which means that all points
$$
z_m-\La_m+s\,,\qquad s=0\lc 2\La_m\,,\quad m=1\lc n\,,
$$
are \pd/.
\Par
Consider the second order \deq/
$$
\tau(u)Q(u)=Q(u-1)\prod_{m=1}^n(u-z_m+\La_m)+
\ka\>Q(u+1)\prod_{m=1}^n(u-z_m-\La_m)\,.
\Tag{Bax}
$$
\wrt/ $Q(u)$. Here $\ka$ is a fixed complex number.
\Par
Let $\Sc_m=\{\,z_m-\La_m+s\vert s=0\lc 2\La_m\,\}$, $m=1\lc n$. Set
$\Sc=\Cup_{m=1}^n\Sc_m$. Let $\F_m=\{\,f:\Sc_m\to\C\,\}$ and let
$\F=\{\,f:\Sc\to\C\,\}$. Let $\pi_m:\F\to\F_m$ be the canonical projection:
$\pi_m\phi=\phi\vst{\Sc_m}\!\!$.
\Par
We consider the next problems related to \deq/ \(Bax).
\demo{\it Global problem, $\ka\ne 1$}
To determine a \pol/ $\tau(u)$ \st/ there exists a nontrivial \pol/ \sol/
to \eq/ \(Bax).
\enddemo
\demo
{\it Global problem, $\ka=1$}
To determine a \pol/ $\tau(u)$ \st/ there exists a nontrivial \pol/ \sol/
to \eq/ \(Bax) of degree at most $1/2+\sum_{m=1}^n\La_m$.
\enddemo
\demo{\it Local problem}
To determine a \pol/ $\tau(u)=(1+\ka)\>u^n+\ldots$  of degree $n$ \st/ there is
$Q\in\F$, $\pi_mQ\ne 0$, $m=1\lc n$, satisfying  \eq/ \(Bax) for all $u\in\Sc$.
\enddemo
\Rem
The specification of the \glob/ for $\ka=1$ is motivated by Lemma 3.3.
\enddemo
\Rem
Note that a given $Q\in\F$ can satisfy \eq/ \(Bax) for at most one \pol/
$\tau(u)=(1+\ka)\>u^n+\ldots$ of degree $n$, since sets $\Sc_m$, $m=1\lc n$,
are \pdj/.
\enddemo
In this section we show (Theorems 3.4, 3.5\>) that \sol/ to the \loc/ is also
a \sol/ to the \glob/ and count the number of local \sol/s.
\Rem
Write equation \(Bax) in the form
$$
\tau(u)Q(u)=\Dlp(u)Q(u-1)+\Dlm(u)Q(u+1)
\Tag{Bax1}
$$
and restrict it to $\Sc_m$.
Then we have a \fd/ homogeneous system of linear \eq/s
$$
\align
\NN4>
\tau(z_m-\La_m)Q(z_m-\La_m)={}& \Dlm(z_m-\La_m)Q(z_m-\La_m+1)\,,
\Tag{Scm}
\\
\tau(z_m-\La_m+s)Q(z_m-\La_m+s)={}& \Dlp(z_m-\La_m+s)Q(z_m-\La_m+s-1)+{}
\\
&{+}\,\Dlm(z_m-\La_m+s) Q(z_m-\La_m+s+1)\,,
\\
\tau(z_m+\La_m)Q(z_m+\La_m)={}&\Dlp(z_m+\La_m) Q(z_m+\La_m-1)\,,
\endalign
$$
$s=1\lc 2\La_m-1$. This system has a nontrivial \sol/ if its determinant is
zero, that is if $\tau\vst{\Sc_m}$ satisfies one additional equation.
\par
Equation \(Bax) motivates the following definition for a general \deq/s of
form \(Bax1): Say that \eq/ \(Bax1) has a \rsp/ at $z$ with
exponents $0,\la$ for $\la\in\Zp$, if $\Dlp(z-\la/2)=\Dlm(z+\la/2)=0$.
It would be interesting if this notion could lead to a \dif/ analog
of the theory of \difl/ \eq/s with \rsp/s.
\enddemo
Now let us return to the local and global problems.
\par
Let $Q\in\F$. Say that $Q$ is a {\it pseudoconstant\/} if all projections
$\pi_mQ$ are constant \fn/s.
\Lm{Q=Q}
Let $\tau(u)$ be a \sol/ to the \loc/. Then a \sol/ $Q\in\F$ to \eq/ \(Bax)
is unique modulo a pseudoconstant factor.
\endpro
\Pf.
It is clear that \eq/ \(Bax) splits into $n$ independent \eq/s for projections
$\pi_1Q\lc\pi_nQ$. The \eq/ for projection $\pi_mQ$ is system \(Scm). Since,
any two of \sol/s to system \(Scm) are proportional, the lemma is proved.
\epf
\Lm{propor}
Let $\tau(u)$ be a \sol/ to the \glob/. Then
\iitem{a)} $\ \deg\tau=n$ and $\tau(u)=(1+\ka)\>u^n +\ldots$.
\iitem{b)} For a given $\tau(u)$, any two of the required \pol/ \sol/s to \eq/
\(Bax) are proportional.
\endpro
\Pf.
Let $\tau(u)=\sum_{k=0}^s a_ku^{s-k}\!$, $a_0\ne 0$. Let
$Q(u)=\sum_{k=0}^\ell b_ku^{\ell-k}\!$, $b_0\ne 0$, be a required \pol/ \sol/
to \eq/ \(Bax). Set $\La=\sum_{m=1}^n\La_m$. %, $\zsi=\sum_{m=1}^n z_m$.
If $\ka\ne 1$, then \eq/ \(Bax) implies claim a) as well as
$$
\align
\ell={}& \La-\tsum_{m=1}^n z_m(1+\ka)/(1-\ka)-a_1/(1-\ka)\,,
\Tag{recu1}
\\
\nn10>
k\>b_k={}& b_0F_k(a_1\lc a_{k+1}, b_1\lc b_{k-1})\,,\qqq k=1\lc\ell\,.
\endalign
$$
If $\ka=1$, then \eq/ \(Bax) implies claim a) as well as
$$
\align
\ell(\ell-2\La-1)={}& a_2+\tsum_{k=2}^n\tsum_{m=1}^{k-1}(z_kz_m+\La_k\La_m)\,,
\Tag{recu2}
\\
\nn10>
k(2\ell-2\La-k-1)\>b_k={}& b_0F_k(a_2\lc a_{k+2}, b_1\lc b_{k-1})\,,\qqq
k=1\lc\ell\,.
\endalign
$$
Both relations \(recu1) and \(recu2) clearly fix $Q(u)$ modulo a constant
factor, since for $\ka=1$ we assume that $\ell\le\La+1/2$.
\epf
Later on, if $\tau(u)$ is a \sol/ to the global (local) problem and $Q(u)$ is
the corresponding \sol/ to \eq/ \(Bax), then we also say that $\tau(u),\,Q(u)$
is a global (local) \sol/.
\Par
Let $\tau(u),\,Q(u)$ be a \sol/ to the \glob/. If $\tau(u)$ is also a \sol/ to
the \loc/ (that is if $Q\vst{\Sc_m}\ne 0$, $m=1\lc n$), then say that
$\tau(u)$ is an {\it \adm/} global \sol/.
\Par
Let $\g=\gg$.
Let $V_1\lc V_n$ be \irr/ \hwm/s with \hw/s $\La_1\lc\La_n$, \resp/.
Set $V=V_1\lox V_n$. Let $\sing$ be the subspace of singular vectors in $V\!$.
\Th{<>1}
Let $\ka$ be generic. Let $z_1\lc z_n$ be \wsep/. Then
\iitem{a)} All \sol/s to the \loc/ are also \sol/s to the \glob/.
\iitem{b)} The number of \sol/s to the \loc/ is equal to $\dim V\!$.
\iitem{c)} If $\tau(u),\,Q(u)$ is an \adm/ \sol/ to the \glob/, then
$\deg Q\le 2\sum_{m=1}^n\La_m$.
\endpro
\Th{=1}
Let $\ka=1$. Let $z_1\lc z_n$ be generic. Then
\iitem{a)} All \sol/s to the \loc/ are also \sol/s to the \glob/.
\iitem{b)} The number of \sol/s to the \loc/ is equal to $\dim V\!$.
\iitem{c)} If $\tau(u),\,Q(u)$ is an \adm/ \sol/ to the \glob/, then
$\deg Q\le \sum_{m=1}^n\La_m$.
\endpro
We prove these theorems in two steps. First we obtain the required number of
\adm/ global \sol/s. Next we  show that the number of local \sol/s cannot
exceed $\dim V$ or $\dim\sing$, \resp/.
\Lm{Basol}
Let $z_1\lc z_n$ be \wsep/. Let $t_1\lc t_\ell$ be an \adm/ \off/ \sol/
to system \(Bae). Let $\tau(u)=\tau(u,t)$ be given by formula \(tau). Set
$Q(u)=\prod_{a=1}^\ell(u-t_a)$. Then $\tau(u),\,Q(u)$ is an \adm/ \sol/ to
the \glob/.
\endpro
\Pf.
By \(tau) $\tau(u)$ is a \raf/ in $u$ with only simple poles at points
$t_1\lc t_\ell$. System \(Bae) means that
$\Res_{u=t_{\rlap{\=$\sss a$}}}\tau(u)=0$, hence
$\tau(u)$ is a \pol/. Equation \(Bax) is fulfilled by the definition of
$\tau(u),\,Q(u)$. Since for $\ka=1$, $\ell>1/2+\sum_{m=1}^n\La_m$ there are no
\adm/ \off/ \sol/s to system \(Bae), we have $\deg Q\le 1/2+\sum_{m=1}^n\La_m$
and hence $\tau(u),\,Q(u)$ is a \sol/ to the \glob/. This global \sol/ is
clearly \adm/ since $Q(u\pm\La_m)\ne 0$, $m=1\lc n$, by the definition of
an \adm/ \sol/ $t_1\lc t_\ell$ to system \(Bae).
\epf
\Lm{number}
Let $M$ be the total number of\/{ \Slorb/s} of \adm/ \sol/s to system
\(Bae) for $\ell=0\lc 2\sum_{m=1}^n\La_m$ altogether.
\iitem{a)} Let $\ka$ be generic. Let $z_1\lc z_n$ be \wsep/. Then $M=\dim V\!$.
Moreover, there are no \adm/ \off/ \sol/s to system \(Bae) for
$\ell>2\sum_{m=1}^n\La_m$.
\iitem{b)} Let $\ka=1$. Let $z_1\lc z_n$ be generic. Then $M=\dim\sing$.
Moreover, there are no \adm/ \off/ \sol/s to system \(Bae) for
$\ell>\sum_{m=1}^n\La_m$.
\endpro
\Pf.
Claim a) follows from Theorem \[ka<>1]. Claim b) follows from Theorems \[ka=1]
and \[details].
\epf
Lemmas \[Basol] and \[number] give the required number of \adm/ \sol/s to the
\glob/. To get an estimate from above for the number of local \sol/s we
consider a spectral problem which can be solved by separation of \var/s.
Equation \(Bax) is the \eq/ for separated \var/s in this problem.
All constructions below are motivated by the functional \Ba/ \Cite{S2}.
\Par
Let $\Fo=\F_1\lox\F_n$. We consider $\Fo$ as a space of \fn/s in $n$ \var/s
$x_1\in\Sc_1\llc x_n\in\Sc_n$. Let $y_k^\pm\in\E(\Fo)$, $k=1\lc n$,
be defined as follows:
$$
\gather
y_k^\pm f(x_1\lc x_n)=f(x_1\lc x_k\mp1\lc x_n)\prod_{m=1}^n(x_k-z_m\pm\La_m)\,.
\\
\Text{Set}
\T(u)=(1+\ka)\prod_{m=1}^n(u-x_m)+\sum_{m=1}^n
\prod_{\tsize{k=1\atop k\ne m}}^n{u-x_k\over x_m-x_k}\cdot
(y_m^+ +\ka\>y_m^-)\,.
\endgather
$$
\Lm*
\back\Cite{S2, Sect.\ 2.4}
Coefficients of the \pol/ $\T(u)$ generate a commutative subalgera in
$\E(\Fo)$.
\endpro
\nt
The proof is straightforward.
\Lm{spectr}
\back\Cite{S2, Sect.\ 2.6}
Let $\tau(u),\,Q(u)$ be a \sol/ to the \loc/. Set
$\Qo=\pi_1Q\lox\pi_nQ\ne 0$. Then $\T(u)\Qo=\tau(u)\Qo$.
Moreover, any \egv/ of $\T(u)$ has the form $\Qo$ for a suitable
\sol/ $\tau(u),\,Q(u)$ to the \loc/.
\endpro
\Pf.
The key observation for the proof of the lemma is that for any $f\in\Fo$
$$
\align
\ifMag\kern2em\fi
\T(x_k)f(x_1\lc x_n)={}&
f(x_1\lc x_k-1\lc x_n)\prod_{m=1}^n(x_k-z_m+\La_m)+{}
\Tagg{Txk}
\\
{}+{}&\ka\>f(x_1\lc x_k+1\lc x_n)\prod_{m=1}^n(x_k-z_m-\La_m)\,.
\endalign
$$
The first claim of the lemma means that for any $(x_1\lc x_n)\in\Sc_1\lx\Sc_n$
$$
\bigl(\T(u)\Qo-\tau(u)\Qo\bigr)\>(x_1\lc x_n)=0\,.
\Tag{Ttau}
$$
Fix $x_1\in\Sc_1\llc x_n\in\Sc_n$. Now \lhs/ above is a \pol/ in $u$ of degree
$n-1$ and it suffices to show that it is zero for $n$ \pd/ values of $u$.
Consider points $x_1\lc x_n$ which are \pd/, since $x_m\in\Sc_m$. Substitution
$u=x_m$ to \(Ttau) reduces it by means of \(Txk) to \eq/ \(Bax) for projection
$\pi_mQ$ at $u=x_m$, which is valid.
\par
Let $f\in\Fo$ be an \egv/ of $\T(u)$ with an \eva/ $\tau(u)$:
$$
\T(u)f(x_1\lc x_n)=\tau(u)f(x_1\lc x_n)\,.
\Tag{Ttauf}
$$
Then $\tau(u)$ is a \pol/ in $u$ of degree $n$ and
$\tau(u)=(1+\ka)\>u^n+\ldots$. Substituting $u=x_1$ to \eq/ \(Ttauf) we get
that $\tau(x_1),\,f(x_1,x_2\lc x_n)$ satisfy \eq/ \(Bax) \wrt/ the \var/
$x_1\in\Sc_1$. Similar to the proof of Lemma \(Q=Q) we obtain that
$f(x_1\lc x_n)=f_1(x_1)\ox f_{[1]}(x_2\lc x_n)$ for $f_1\in\F_1$ which obeys
\eq/ \(Bax) and suitable $f_{[1]}\in\F_2\lox\F_n$.
Similarly, $\tau(x_2),\,f_{[1]}(x_2,x_3\lc x_n)$ satisfy \eq/ \(Bax) \wrt/
a \var/ $x_2\in \Sc_2$ and
$f_{[1]}(x_2\lc x_n)=f_2(x_2)\ox f_{[2]}(x_3\lc x_n)$ for
suitable $f_2,\,f_{[2]}$, etc. Finally, $f(x_1\lc x_n)=f_1(x_1)\lox f_n(x_n)$
and $\tau(u),\,f_m(u)$ satisfy \eq/ \(Bax) for $u\in\Sc_m$, $m=1\lc n$. Since
$\Sc_k\cap\Sc_m=\Empty$ unless $k=m$, there is a \sol/ $\tau(u),\,Q(u)$
to the \loc/, \st/ $f_1=\pi_1Q\llc f_n=\pi_nQ$ which means $f=\Qo$.
\epf
\Pf of Theorem \[<>1].
By Lemmas \[Basol], \[number] and \[propor] we point out $\dim V\!$ \adm/
\sol/s $\tau(u),\,Q(u)$ to the \glob/. For all of them
$\deg Q(u)\le 2\sum_{m=1}^n\La_m$.
Since any \adm/ global \sol/ is also a local \sol/, Lemma \[spectr] shows that
we exhaust all local as well as \adm/ global \sol/s.
\epf
 From now until the end of the section, let $\ka=1$. To get the required
estimate from above for the number of \sol/s to the \loc/ in this case
we equip the space $\Fo$ with a structure of a \gmod/ isomorphic to the
\gmod/ $V\!$. Set
$$
\align
H &{}=\sum_{m=1}^n\Bigl(x_m-z_m-\prod_{\tsize{k=1\atop k\ne m}}^n
(x_m-x_k)\1\cdot y_m^-\Bigr)\,,
\\
F &{}= \sum_{m=1}^n\prod_{\tsize{k=1\atop k\ne m}}^n(x_m-x_k)\1\cdot y_m^-\,,
\\
E &{}=\sum_{m=1}^n\Bigl(2(x_m-z_m)-\prod_{\tsize{k=1\atop k\ne m}}^n(x_m-x_k)\1
\cdot(y_m^+ + y_m^-)\Bigr)\,,
\\
\Dl(u) &{}=\prod_{m=1}^n(u-z_m+\La_m)(u-z_m-\La_m-1)\,.
\endalign
$$
\Lm{Fmod}
The map $e\to E$, $f\to F$, $h\to H$ makes the space $\Fo$ into a \gmod/
isomorphic to the \gmod/ $V\!$.
\endpro
\Pf.
Verifying commutation relations between $E,\,F,\,H$ is cumbersome, but
straightforward. The following identity is useful:
$$
\align
\sum_{m=1}^n\Bigl( \Dl(x_m+1)&
\prod_{\tsize{k=1\atop k\ne m}}^n{1\over(x_m-x_k)(x_m-x_k+1)}\,-{}
\\
\nn3>
{}-\Dl(x_m) &\prod_{\tsize{k=1\atop k\ne m}}^n
{1\over(x_m-x_k-1)(x_m-x_k)}\,\Bigr)\,=\,
2\sum_{m=1}^n(x_m-z_m)\,.
\endalign
$$
The identity itself is \eqv/ to
$$
\Bigl(\,\Res_{u=\8}+\sum_{m=1}^n(\>\Res_{\,u=x_m}+\Res_{u=x_m+1}\>)\Bigr)\,\,
\Dl(u)\prod_{m=1}^n{1\over(u-x_m)(u-x_m-1)}\,=\,0\,.
$$
The character of the obtained \gmod/ $\Fo$ coincides with the character of the
\gmod/ $V$ since dimensions of the corresponding \wt/ subspaces are clearly the
same. This proves the required isomorphism.
\epf
Define one more \pol/ taking values in $\E(\Fo)$:
$$
\align
& \aligned
\Cc(u)=\sum_{m=1}^n\biggl(\,
\prod_{\tsize{k=1\atop k\ne m}}^n{u-x_k\over x_m-x_k}
\cdot\Bigl(\Dl(x_m+1)&\prod_{k=1}^n(x_m-x_k+1)\1 +{}
\\
{}+\Dl(x_m)&\prod_{k=1}^n(x_m-x_k-1)\1 - y_m^+ - y_m^-\Bigr)+{}
\endaligned
\\
\nn2>
\ald
&\hp{C(u)=\sum_{m=1}^n\biggl(\,}\llap{$+$}
\sum_{\tsize{l=1\atop l\ne m}}^n (x_m-x_l)\1(x_m-x_l-1)\1
\prod_{\tsize{k=1\atop k\ne l,m}}^n{u-x_k\over(x_m-x_k)(x_l-x_k)}\cdot
y_m^+\>y_l^-\biggr)\,.
\endalign
$$
The next lemma is proved by tremendous, but all the same straightforward
calculation.
\Lm{commute}
\iitem{a)} Coefficients of the \pol/ $\T(u)+\Cc(u)$ generate a commutative
subalgera in $\E(\Fo)$.
\iitem{b)} Coefficients of the \pol/ $\T(u)+\Cc(u)$ commute with the
$\gg$ action in $\Fo$.
\iitem{c)}  Coefficients of the \pol/ $\Cc(u)$ are raising operators:
$[H,\Cc(u)]=H$.
\endpro
\Pf of Theorem \[=1].
By Lemmas \[Fmod] and \[commute] $\T(u)+\Cc(u)$ can have at most $\dim\sing$
different \eva/s and the same is true for $\T(u)$. Hence, by Lemma \[spectr]
there are at most $\dim\sing$ \sol/s to the \loc/.
\par
On the other hand, by Lemmas \[Basol], \[number] and \[propor] we point out
$\dim\sing$ \adm/ \sol/s $\tau(u),\,Q(u)$ to the \glob/. For all of them
$\deg Q(u)\le \sum_{m=1}^n\La_m$. Since any \adm/ global \sol/ is also
a local \sol/, we exhaust all local as well as \adm/ global \sol/s.
\epf

\Sect{\=Bases of \Bv/s in $\UU$-modules}
In this section we describe a \^{$q$-variant} of Theorem \[ka<>1].
Proofs of these new statements are completely similar to the corresponding
proofs of Section 1. Notations, used in this section differ slightly from
the notations, used in Section 1.
\Par
Let $\Co=\C\setminus\{0\}$. Let $q\in\Co$, $q^2\ne 1$.
Let $\g=\gg$. Let $e,f,q^h,q^{-h}$ be generators of $\Uq$:
$$
q^he\>\>\!q^{-h}=q\>\>\!e\,,\qqq q^hfq^{-h}=q\1f\,,\qqq
[e,f]={q^{2h}-q^{-2h}\over q-\q}\,.
$$
Let $\M=\E(\C^2)$. Introduce $T(u)\in\M[u]\ox\Uq$ as follows:
$$
T(u)=\pmatrix uq^h-q^{-h}& uf\>(q-\q)\\
\nn4> e\>(q-\q)& uq^{-h}-q^h\endpmatrix\,.
$$
Let $\io_m$ be the embedding $\Uq\to\Uqn$ as the \^{$m$-th} tensor factor.
\par
Let $z=(z_1\lc z_n)\in\Con$.
Let $T_m(u)=\id\ox\io_m\bigl(T(u)\bigr)\in\M[u]\ox\Uqn$. Set
$$
\TT(u)=
z_1T_1(u/z_1)\ldots z_nT_n(u/z_n)=\pmatrix A(u)& B(u)\\ C(u)& D(u)\endpmatrix
$$
where $A(u)$, $B(u)$, $C(u)$ and $D(u)$ are suitable elements in $\Cu\ox\Un$.
\par
Let $\ka\in\C$. Set $\T(u)=A(u)+\ka D(u)$. Coefficients of the \pol/ $\T(u)$
generate a remarkable commutative subalgebra in $\Un$.
\Rem
 From the physical point of view the \^{$q$-deformed} case considered in
this section corresponds to partially anisotropic quantum models, while
the previous case corresponds to isotropic quantum models.
\enddemo
Later on in this section we assume that $q$ is not a root of unity.
In the next section we partially extend results of this section to the case
when $q$ is a root of unity.
\Par
Let $V_1\lc V_n$ be \irr/ \hwu/s with \hw/s $\La_1\lc\La_n$ and \gv/s
$v_1\lc v_n$, \resp/. Set $V=V_1\lox V_n$. Set $\La=\sum_{m=1}^n\La_m$.
Let $\ell\in\Zp$. Let $\Vl\sub V$ be a \wt/ subspace
$\Vl=\{\,v\in V\vert q^h\lox q^h\cdot v=q^{\La-\ell}v\,\}$.
\Rem
All the time only $q^{\La_1}\lc q^{\La_n}\!$ are used. $\La_1\lc \La_n$
themselves never appear in formulae.
\enddemo
Let $t=(t_1\lc t_\ell)\in\Col$. Set
$$
w(t)=B(t_1)\ldots B(t_\ell)\cdot v_1\lox v_n\,.
\Tag{BBB}
$$
$w(t)$ is a \Vlv/ \sym/ \pol/ in \var/s $t_1\lc t_\ell$. Another formula for
$w(t)$ see below in Lemma 4.10. %% !!! %%
\par
For given $z_1\lc z_n$ consider a system of
algebraic equations on \var/s $t_1\lc t_\ell$:
\ifMag\vv->\fi
$$
\prod_{m=1}^n(q^{2\La_m}t_a-z_m)\,
\prod_{\tsize{b=1\atop b\ne a}}^\ell(t_a-q^2t_b)=
\ka\prod_{m=1}^n(t_a-q^{2\La_m}z_m)\,
\prod_{\tsize{b=1\atop b\ne a}}^\ell(q^2t_a-t_b)\,.
\Tag{Baeq}
$$
$a=1\lc\ell$. We consider this system only in $\Col$.
A \sol/ $t$ to system \(Baeq) is called an {\it\off/\/} \sol/ if
are \pd/, and a {\it diagonal\/} \sol/, otherwise.
\par
Set\vv->
\ifMag
$$
\align
\tau(u,t)=q^{-\La-\ell} &\Bigl(\>\prod_{m=1}^n(q^{2\La_m}u-z_m)\,
\prod_{a=1}^\ell{u-q^2t_a\over u-t_a}+{}
\Tag{tauq}
\\
&{}+\ka\prod_{m=1}^n(u-q^{2\La_m}z_m)\,
\prod_{a=1}^\ell{q^2u-t_a\over u-t_a}\>\Bigr)\,.
\endalign
$$
\else
$$
\tau(u,t)=q^{-\La-\ell}\Bigl(\>\prod_{m=1}^n(q^{2\La_m}u-z_m)\,
\prod_{a=1}^\ell{u-q^2t_a\over u-t_a}+
\ka\prod_{m=1}^n(u-q^{2\La_m}z_m)\,
\prod_{a=1}^\ell{q^2u-t_a\over u-t_a}\>\Bigr)\,.
\Tag{tauq}
$$
\fi
\Th{eigenq}
\back\Cite{BIK},\,\Cite{F},\,\Cite{FT}
Let $t_1\lc t_\ell$ be an \off/ \sol/ to system \(Baeq).
\ifMag\else\nl\fi
Then $\T(u)\cdot w(t)=\tau(u,t)\>w(t)$.
\endpro
Define a set $\ZZ$ by the \eq/
$$
\prod_{a=1}^\ell\,\Bigl(\,\prod_{m=1}^n
(q^{2\La_m}t_a-z_m)(t_a-q^{2\La_m}z_m)
\prod_{\tsize{b=1\atop b\ne a}}^\ell(t_a-q^2t_b)\Bigr)=0\,.
$$
A \sol/ $t$ to system \(Baeq) is called an {\it \adm/\/} \sol/ if $t\nin\ZZ$
and an {\it un\adm/\/} otherwise. For an \adm/ \off/ \sol/ $t$ the vector
$w(t)$ is called the {\it \Bv/\/}. A \sol/ $t$ to system \(Baeq) is called a
{\it trivial\/} \sol/ if $w(t)=0$ and {\it nontrivial\/} otherwise.
\par
System \(Baeq) is preserved by the natural action of the \symg/ $\Sl$ on \var/s
$t_1\lc t_\ell$. Therefore, $\Sl$ acts on \sol/s to this system.
Let $\CC$ be the set of \Slorb/s of \adm/ \off/ \sol/s.
\par
Say that $z_1\lc z_n$ are {\it \wsep/\/} if all points
$z_mq^{2(s-\La_m)}\!$, $s\in\Zp$, $s\le\dim V_m-2$, and $z_mq^{2\La_m}\!$,
$m=1\lc n$, are \pd/.
\Rem
$z_1\lc z_n$ are \wsep/ if and only if module $V$ enjoys next properties:
\iitem{i)} $V$ is \irr/ \wrt/ a subalgebra \gb/ coefficients of \pol/s
$A(u),\alb\,B(u),\alb\,C(u),\alb\,D(u)$.
\iitem{ii)} A commutative subalgebra \gb/ coefficients of $A(u)$ acts in $V$
in a semisimple way.
\par\nt
This follows from results of \Cite{T},\,\Cite{NT}.
The part ``{\sl only if\/}'' also follows from results of the present paper.
Note in addition that the action in $V$
of subalgebra \gb/ coefficients of $A(u)$ has simple spectrum.
\enddemo
\Th{kappaq}
Let $\ka$ be generic. Let $z_1\lc z_n$ be \wsep/. Then
\iitem{a)} All \adm/ \off/ \sol/s to system \(Baeq) are \ndeg/.
\iitem{b)} All un\adm/ \off/ \sol/s to system \(Baeq) are trivial.
\iitem{c)} $\#\CC=\dim\Vl$ and the corresponding \Bv/s form a base in $\Vl$.
\endpro
\nt
The proof is completely similar to the proof of Theorem \[ka<>1].
We give below only the main points of the proof.
\Rem
This Theorem was proved in \Cite{TV} for $q$ not a root of unity and generic
$z_1\lc z_n$, $\La_1\lc\La_n$. The proof of Theorem \[ka<>1] is similar to
the corresponding proof in \Cite{TV}.
\enddemo
\Rem
A theorem similar to Theorem \[kappaq] for the case $\La_m=1/2$, $m=1\lc n$,
was announced in a recent preprint \Cite{LS}.
\enddemo
\Lm{admiq}
Let $q^{4\La_m}\nin\{\,q^{2s}\vert s\in\Zp\,\}$, $m=1\lc n$. For $\ka\ne 0$
and \wsep/ $z_1\lc z_n$ all \sol/s to system \(Bae) are \adm/.
\endpro
\Lm{isolq}
System \(Baeq) has no nonisolated \off/ \sol/s unless
$\ka=q^{2(s-\ell+\La)}\alb=q^{2(\ell-\s-\La)}$ for some
$s,\s\in\{\,1\lc\ell\,\}$.
\endpro
\Cr*
System \(Baeq) has no nonisolated \off/ \sol/s unless
$q^{4\La}\in\{\,1,q^2\lc q^{4\ell-4}\,\}$.
\endpro
\Lm{88offq}
Consider \sol/s to system \(Baeq) as (multivalued) \fn/s of $\ka$.
\iitem{a)} For any \off/ \sol/ to system \(Baeq) every its branch remains
finite for any $\ka\ne q^{2(s-\ell+\La)}$, $s=1\lc\ell$.
\iitem{b)} For any \off/ \sol/ to system \(Baeq) every its branch remains
in $\Col$ for any $\ka\ne q^{2(\ell-s-\La)}$, $s=1\lc\ell$.
\endpro
\Lm{offq}
\iitem{a)} All \sol/s to system \(Baeq) at $\ka=0$ are isolated.
\iitem{b)} For generic $z_1\lc z_n$ all \off/ \sol/s to system \(Baeq) at
$\ka=0$ are \ndeg/.
\endpro
\Lm{ddq}
Let $z_1\lc z_n$ be generic.
Let $t(\ka)$ be a \sol/ to system \(Baeq), which is a deformation of
a diagonal \sol/ $t(0)$ to system \(Baeq) at $\ka=0$. Then $t(\ka)$ is
a diagonal \sol/.
\endpro
Set $\Zl=\{\,\nu\in\Z^n_{\ge 0}\vert\sum_{m=1}^n\nu_m=\ell\,\}$.
\vv.1>
Set $\Zlo=\{\,\nu\in\Z^n_{\ge 0}\vert\sum_{m=1}^n\nu_m\alb=\ell$ and
$\nu_m<\dim V_m$, $m=1\lc n\,\}$. For $\nu\in\Zl$ let $\ts(\nu)$ be the
following \off/ \sol/ to system \(Baeq) at $\ka=0$
$$
\ts_a(\nu)=q^{2(a-\ell_{i-1}-1-\La_i)}z_i\qqq\for \ell_{i-1}<a\le\ell_i
\Tag{tnuq}
$$
where $\ell_i=\sum_{m=1}^i\nu_m$, \,$\ell_0=0$, $\ell_n=\ell$.
Let $t(\nu,\ka)$ be a \sol/ to system \(Baeq) which is a deformation of
$\ts(\nu)$.
\Lm{Zloq}
Let $\ka$ be generic.
\iitem{a)} Let $\nu\in\Zlo$. Let $z_1\lc z_n$ be \wsep/.
Then $t(\nu,\ka)$ is an \adm/ \sol/.
\iitem{b)} Let $\nu\in\Zl\setminus\Zlo$. Let $z_1\lc z_n$ be generic.
Then $t(\nu,\ka)$ is an un\adm/ \sol/.
\endpro
Introduce the canonical monomial base in $V$:
$\{\,F^\nu = f^{\nu_1}v_1\lox f^{\nu_n}v_n\,\}$. It is clear that
$\{\,F^\nu\vert\nu\in\Zlo\,\}$ is a base in $\Vl$.
\Lm{decompq}
\back\Cite{Ko}
The following decomposition holds
\ifMag
$$
\align
w(t)=\sum_{\text{\it part}}\,F^{\nu(\Gm)}\,&\prod_{l=2}^n\,\prod_{m=1}^{l-1}\,
\Bigl(
\prod_{\tsize{a\in\Gm_l\atop\rlap{\=$\ssize b\in\Gm_m$}\hp{a\in\Gm_l}}}
{q\1t_a-qt_b\over t_a-t_b}\x{}
\\
{}\x{}&\prod_{\slap{\=$\ssize a\in\Gm_l$}}^{\hp{l=2}}
\,(q^{\La_m}t_a-q^{-\La_m}z_m)
\prod_{a\in\Gm_m}(q^{-\La_m}t_a-q^{\La_m}z_m)\Bigr)\,.
\endalign
$$
\else
$$
w(t)=\sum_{\text{\it part}}\,F^{\nu(\Gm)}\,\prod_{l=2}^n\,\prod_{m=1}^{l-1}\,
\Bigl(
\prod_{\tsize{a\in\Gm_l\atop\rlap{\=$\ssize b\in\Gm_m$}\hp{a\in\Gm_l}}}
{q\1t_a-qt_b\over t_a-t_b}
\prod_{a\in\Gm_l}(q^{\La_m}t_a-q^{-\La_m}z_m)
\prod_{a\in\Gm_m}(q^{-\La_m}t_a-q^{\La_m}z_m)\Bigr)\,.
$$
\fi
Here the sum is taken over all partitions of the set $\{1\lc\ell\}$ into
disjoint subsets $\Gm_1\lc\Gm_n$ and $\nu(\Gm)=(\#\Gm_1\lc\#\Gm_n)$.
\endpro
Let $\v_m$ be a linear \fn/ on $V_m$ \st/ $\bra\v_m,v_m\ket=1$ and
$\bra\v_m,v\ket=0$ for any \wt/ vector $v\in V_m$, $v\ne v_m$.
\Th{orthoq}
\back\Cite{Ko},\,\Cite{TV}
Let $\La_1\lc\La_n$, $z_1\lc z_n$ and $\ka$ be generic.
\iitem{a)} For any \off/ \sol/ $t=(t_1\lc t_n)$ to system \(Baeq)
\par
\ifMag
$$
\align
\bra\v_1\lox{}&\v_n,C(t_1)\ldots C(t_n)w(t)\ket={}
\\
\nn4>
&\aligned
\llap{${}={}$}\>(-1)^\ell q^{-2\La-\ell(\ell+1)}(q-\q)^\ell
\prod_{m=1}^n\,\prod_{a=1}^\ell(t_a-q^{2\La_m}z_m)\,
\prod_{a=2}^\ell\,\prod_{b=1}^{a-1}{q^2t_a-t_b\over t_a-t_b}\>\x{}&
\\
\nn6>
{}\x\>\det\left[
\aligned
t_a\ddt\Bigl(\> &\prod_{m=1}^n(q^{2\La_m}t_b-z_m)\,
\prod_{\tsize{c=1\atop c\ne b}}^\ell(t_b-q^2t_c)-{}
\\
{}-\ka &\prod_{m=1}^n(t_b-q^{2\La_m}z_m)\,
\prod_{\tsize{c=1\atop c\ne b}}^\ell(q^2t_b-t_c)\Bigr)
\endaligned
\right]_{a,b=1\lc\ell}\,& .
\endaligned
\endalign
$$
\else
$$
\align
\bra\v_1 &{}\lox\v_n,C(t_1)\ldots C(t_n)w(t)\ket={}
\\
\nn4>
&{}=(-1)^\ell q^{-2\La-\ell(\ell+1)}(q-\q)^\ell
\prod_{m=1}^n\,\prod_{a=1}^\ell(t_a-q^{2\La_m}z_m)\,
\prod_{a=2}^\ell\,\prod_{b=1}^{a-1}{q^2t_a-t_b\over t_a-t_b}\,\x{}
\\
\nn6>
&\>{}\x\>\det\Biggl[t_a\ddt\Bigl(\>\prod_{m=1}^n(q^{2\La_m}t_b-z_m)\,
\prod_{\tsize{c=1\atop c\ne b}}^\ell(t_b-q^2t_c)-
\ka\prod_{m=1}^n(t_b-q^{2\La_m}z_m)\,
\prod_{\tsize{c=1\atop c\ne b}}^\ell(q^2t_b-t_c)\Bigr)
\Biggr]_{a,b=1\lc\ell}\,.
\endalign
$$
\fi
\iitem{b)} For any \off/ \sol/s $t$ and $\t=(\t_1\lc\t_n)$ to system \(Bae)
which lie on different \Slorb/s
$$
\bra\v_1\lox\v_n,C(\t_1)\ldots C(\t_n)w(t)\ket=0\,.
$$
\iitem{c)} Let $t,\t$ be isolated \sol/ to system \(Bae). Then both claims
{\rm a)} and {\rm b)} remain valid for any $\La_1\lc\La_n$, $z_1\lc z_n$ and
$\ka$.
\endpro
\Rem
In this paper we use a normalization of $w(t)$ which differs from the
normalization in \Cite{TV}.
\enddemo
\Lm{sigmaq}
\back{\Cite{T} {\rm (cf.\ \Cite{CP2} for more detailed proof)}}
Let $\si$ be a \perm/ of $1\lc n$. Set $V^\si=V_{\si(1)}\lox V_{\si(n)}$
and $z^\si=(z_{\si(1)}\lc z_{\si(n)})$. Let $w^\si(t)$ be constructed in the
same manner for $V^\si, z^\si$ as $w(t)$ is constructed for $V,z$. Then
there is a linear isomorphism $V\to V^\si$ \st/ $w(t)\map w^\si(t)$.
\endpro

\Sect{\=Bases of \Bv/s in $\UU$-modules at roots of unity}
Let $q$ be a root of unity. In this case a \rep/ theory of $\Uq$ drastically
changes (see e.g.\ \Cite{DCK}\>). Nevertheless, results of the previous section
essentially remain valid in this case. We give precise statements in this
section.
\par
Later on we keep notations from the previous section unless changes are given
explicitly.
\Par
Let $q^2$ be a primitive \^{$N$-th} root of unity. Let V be an \irr/ \hwm/
with \hw/ $\La$. In the case in question there are two types of such modules.
\par\hglue.5\parindent
\llap{I.\enspace}{\it Restricted modules.} These modules correspond to
$q^{4\La}\in\{\,1,\,q^2,\lc q^{2N-4}\}$. A restricted module $V$ is uniquely
fixed by its \hw/ $\La$;  $q^{2\dim V}=q^{4\La+2}$, $\dim V<N$. All these
modules admit a deformation to the case of generic $q$.
\par\hglue.5\parindent
\llap{II.\enspace}{\it Unrestricted modules.} These modules correspond to
$q^{4\La}\in\Co\setminus\{\,1,\,q^2,\lc q^{2N-4}\}$. In this case $\dim V=N$.
For a given \hw/ $\La$ there is a one-parametric family of unrestricted
modules. They are separated by values of the central element $f^N$ in these
modules. The value of $f^N$ in an unrestricted module can be an arbitrary
complex number. The only unrestricted modules, which can be deformed to the
case of generic $q$, are those corresponding to $q^{4\La}=q^{-2}\!$, $f^N=0$.
We call these modules {\it quasirestricted}.
\Par
Let $V_1\lc V_n$ be \irr/ \hwu/s with \hw/s $\La_1\lc\La_n$ and \gv/s
$v_1\lc v_n$, \resp/. Set $V=V_1\lox V_n$. Set $\La=\sum_{m=1}^n\La_m$.
Let $\ell\in\Zp$. Let $\Vl\sub V$ be the subspace  \gb/
$\{\,f^{\nu_1}v_1\lox f^{\nu_n}v^n\vert\sum_{m=1}^n\nu_m=\ell\,\}$.
\par
Let $t=(t_1\lc t_\ell)\in\Col$. Let $w(t)$ be defined by \(BBB). It is a \Vlv/
\sym/ \pol/ in \var/s $t_1\lc t_\ell$. Let $\tau(u,t)$ be defined by \(tauq).
Theorem \[eigenq] is known to remain valid even if $q$ is a root of unity.
\par
Let $\CC$ be the set of \Slorb/s of \adm/ \off/ \sol/s to system \(Baeq).
\par
Set $s_m=\dim V_m-2$ if $V_m$ is a restricted or quasirestricted module and
$s_m=N-1$, otherwise. Say that $z_1\lc z_n$ are {\it \wsep/\/} if all points
$z_mq^{2(s-\La_m)}\!$, $s=0\lc s_m$, and $z_mq^{2\La_m}\!$, $m=1\lc n$,
are \pd/.
\Th{kappar}
Let $\ka$ be generic. Let $z_1\lc z_n$ be \wsep/. Then
\iitem{a)} All \adm/ \off/ \sol/s to system \(Baeq) are \ndeg/.
\iitem{b)} $\#\CC=\dim\Vl$ and the corresponding \Bv/s form a base in $\Vl$.
\endpro
\Lm{admisol}
Let $\ka\ne 0$. Then all \adm/ \sol/s to system \(Baeq) are isolated provided
the couple $\ka,q^{\La}$ does not belong to a certain finite set.
\endpro
\Pf.
Assume that there is a nonisolated \adm/ \sol/ to system \(Baeq). This means
that we have a curve $t(s)$, $s\in\R$, such that $t(s)$ is an \adm/ \sol/ to
system \(Baeq) for any $s$. Moreover, we can assume that as $s\to+\8$, $t(s)$
tends to infinity in the following way: $t_a(s)\to\8$ if $a\le f$ and $t_a(s)$
has a finite limit if $a>f$. Taking the product of the first $f$ \eq/s of
system \(Baeq) we obtain
$$
\align
& \prod_{a=1}^f\prod_{\tsize{b=1\atop b\ne a}}^\ell
\bigl(t_a(s)-q^2t_b(s)\bigr)\,
\prod_{a=1}^f\,\Bigl(\,\prod_{m=1}^n\bigl(q^{2\La_m}t_a(s)-z_m\bigr)
\prod_{b=f+1}^\ell\bigl(t_a(s)-q^2t_b(s)\bigr)\Bigr)=
\\
& \ =\>\ka^f\prod_{a=1}^f\prod_{\tsize{b=1\atop b\ne a}}^\ell
\bigl(q^2t_a(s)-t_b(s)\bigr)\,
\prod_{a=1}^f\,\Bigl(\,\prod_{m=1}^n\bigl(t_a(s)-q^{2\La_m}z_m\bigr)
\prod_{b=f+1}^\ell\bigl(q^2t_a(s)-t_b(s)\bigr)\Bigr)\,.
\endalign
$$
The first products in the left and right hand sides above coincide. Moreover,
they are not zero, since $t(s)$ is an \adm/ \sol/. Cancelling these products
and taking the limit $s\to\8$ we obtain that
$$
\ka^f=q^{2f(f-\ell+\La)}\,.
\Tag{kaLa1}
$$
System \(Baeq) is invariant under the transformation $z_m\to z_m\1$,
$t_a\to t_a\1$, $\ka\to\ka\1$, $m=1\lc n$, $a=1\lc\ell$. Therefore we also
have that
$$
\ka^{\f}=q^{2\f(\f-\ell+\La)}
\Tag{kaLa2}
$$
for a suitable $\f\in\{1\lc\ell\}$. Hence, if the couple $\ka,q^\La$ does not
obey \eq/s \(kaLa1) and \(kaLa2) for some $f,\f\in\{1\lc\ell\}$ then all
\adm/ \sol/s to system \(Baeq) are isolated.
\epf
\Pf of Theorem \[kappar].
Let $\ka$ be generic. Then $\ts(\nu)$, $\nu\in\Zlo$, are \ndeg/ \off/ \sol/s
to system \(Baeq) at $\ka=0$. Their deformations $t(\nu,\ka)$, $\nu\in\Zlo$,
are \adm/ \ndeg/ \off/ \sol/s to system \(Baeq) and the corresponding \Bv/s
form a base in $\Vl$. The proof is the same as for generic $q$.
\par
Any unrestricted module can be considered as a continuous deformation of a
quasirestricted module. Since for generic $\ka$ all \adm/ \sol/s are isolated,
it suffices to proof the theorem for the case when all \gmod/s $V_1\lc V_n$
are restricted or quasirestricted. In the last case, all modules $V_1\lc V_n$
can be deformed to the case of generic $q$. Then inequality $\#\CC\le\dim\Vl$
follows from Theorem \[kappaq], which completes the proof of Theorem \[kappar].
\epf

\Sect{Difference \eq/s with \rsp/s. Multiplicative case}
In this section we describe a \^{$q$-variant} of Theorem \[<>1].
Proofs of these new statements are completely similar to the corresponding
proofs of Section 3.  So we give only the most important points of the proofs.
We keep notations used in the two previous sections which differ slightly from
the notations, used in Section 3.
\par
Let $q\in\Co$, $q^2\ne 1$.
Details of consideration depend on whether $q$ is or is not a root of unity.
We give necessary specifications for the case of $q$ being a root of unity
at the end of the section.
\Par
Let $q$ be not a root of unity. In this case we assume that $\La_1\lc \La_n$
are \st/ $q^{4\La_m}=q^{2d_m}$, for some numbers $d_m\in\Zp$, $m=1\lc n$.
Note that the integers $d_1\lc d_n$ are uniquely determined. We also assume
that $z_1\lc z_n$ are \wsep/, which means that all points
$z_mq^{2(s-\La_m)}\!$, $s=0\lc d_m-1$, and $z_mq^{2\La_m}$, $m=1\lc n$,
are \pd/.
\Par
Consider the second order \deq/
$$
\tau(u)Q(u)=Q(q^{-2}u)\prod_{m=1}^n(q^{2\La_m}u-z_m)+
\thi\>Q(q^2u)\prod_{m=1}^n(u-q^{2\La_m}z_m)\,.
\Tag{Baxq}
$$
\wrt/ $Q(u)$. Here $\thi$ is a fixed complex number.
\Par
Let $\Sc_m=\{\,z_mq^{2(s-\La_m)}\vert s=0\lc d_m\,\}$, $m=1\lc n$. Set
$\Sc=\Cup_{m=1}^n\Sc_m$. Let $\F_m=\{\,f:\Sc_m\to\C\,\}$ and let
$\F=\{\,f:\Sc\to\C\,\}$. Let $\pi_m:\F\to\F_m$ be the canonical projection:
$\pi_m\phi=\phi\vst{\Sc_m}\!\!$. Set $\La=\sum_{m=1}^n\La_m$.
\Par
We consider the next problems related to \deq/ \(Baxq).
\demo{\it Global problem}
To determine a \pol/ $\tau(u)$ \st/ there exists a nontrivial \pol/ \sol/
$Q(u)$ to \eq/ \(Baxq) \st/ $Q(0)\ne 0$.
\enddemo
\demo{\it Local problem}
To determine a \pol/ $\tau(u)$ of degree at most $n$,
$\tau(0)=(-1)^n(1+q^{2\La}\thi)\prod_{m=1}^nz_m$, \st/ there is $Q\in\F$,
$\pi_mQ\ne 0$, $m=1\lc n$, satisfying  \eq/ \(Baxq) for all $u\in\Sc$.
\enddemo
\Rem
Note that a given $Q\in\F$ can satisfy \eq/ \(Baxq) for at most one \pol/
$\tau(u)$ of degree $n$ with the prescribed value $\tau(0)$, since the sets
$\Sc_m$, $m=1\lc n$, are \pdj/.
\enddemo
Let $Q\in\F$. Say that $Q$ is a {\it pseudoconstant\/} if all projections
$\pi_mQ$ are constant \fn/s.
\Lm{Q=Qq}
Let $\tau(u)$ be a \sol/ to the \loc/. Then a \sol/ $Q\in\F$ to \eq/ \(Baxq)
is unique modulo a pseudoconstant factor.
\endpro
\Lm{proporq}
Let $\tau(u)$ be a \sol/ to the \glob/. Then
\iitem{a)} $\ \deg\tau\le n$ and
$\tau(0)=(-1)^n(1+q^{2\La}\thi)\prod_{m=1}^nz_m$.
\iitem{b)} Let $\thi\ne q^{-2(s+\La)}$, $s\in\Zpp$. Then for a given $\tau(u)$,
any two of the required \pol/ \sol/s to \eq/ \(Baxq) are proportional.
\endpro
Let $\tau(u),\,Q(u)$ be a \sol/ to the \glob/. If $\tau(u)$ is also a \sol/ to
the \loc/ (that is if $Q\vst{\Sc_m}\ne 0$, $m=1\lc n$), then say that
$\tau(u)$ is an {\it \adm/} global \sol/.
\Th{<>q}
Let $\thi$ be generic. Let $z_1\lc z_n$ be \wsep/. Then
\iitem{a)} All \sol/s to the \loc/ are also \sol/s to the \glob/.
\iitem{b)} The number of \sol/s to the \loc/ is equal to
$\prod_{m=1}^n(d_m+1)$.
\iitem{c)} If $\tau(u),\,Q(u)$ is an \adm/ \sol/ to the \glob/, then
$\deg Q\le \sum_{m=1}^nd_m$.
\endpro
\Lm{Basolq}
Let $z_1\lc z_n$ be \wsep/. Let $t_1\lc t_\ell$ be an \adm/ \off/ \sol/
to system \(Baeq) with $\ka=q^{2\ell}\thi$. Let
$\tau(u)=q^{\La-\ell}\tau(u,t)$ where $\tau(u,t)$ is given by formula \(tauq).
Set $Q(u)=\prod_{a=1}^\ell(u-t_a)$. Then $\tau(u),\,Q(u)$ is an \adm/ \sol/ to
the \glob/.
\endpro
\Lm{numberq}
Let $M$ be the total number of\/{ \Slorb/s} of \adm/ \sol/s to system
\(Baeq) for $\ell=0\lc \sum_{m=1}^n d_m$ altogether. Let $\ka$ be generic.
Let $z_1\lc z_n$ be \wsep/. Then $M=\prod_{m=1}^n(d_m+1)$. Moreover, there are
no \adm/ \off/ \sol/s to system \(Baeq) for $\ell>\sum_{m=1}^nd_m$.
\endpro
\Rem
Since for generic $\ka$ the number of nontrivial \sol/s to system \(Baeq)
does not depend on $\ka$, the explicit dependence $\ka=q^{2\ell}\thi$ is
irrelevant.
\enddemo
Lemmas \[Basolq] and \[numberq] give the required number of \adm/ \sol/s to
the \glob/. To get an estimate from above for the number of local \sol/s we
consider a spectral problem which can be solved by separation of \var/s.
Equation \(Baxq) is the \eq/ for separated \var/s in this problem.
\Par
Let $\Fo=\F_1\lox\F_n$. We consider $\Fo$ as a space of \fn/s in $n$ \var/s
$x_1\in\Sc_1\llc x_n\in\Sc_n$. Let $y_k^\pm\in\E(\Fo)$, $k=1\lc n$,
be defined as follows:
$$
\gather
{\align
y_k^+ f(x_1\lc x_n)  &{}=f(x_1\lc q^{-2}x_k\lc x_n)
\prod_{m=1}^n(q^{2\La_m}x_k-z_m)\,,
\\
y_k^- f(x_1\lc x_n)  &{}=f(x_1\lc q^2x_k\lc x_n)
\prod_{m=1}^n(x_k-q^{2\La_m}z_m)\,.
\endalign}
\\
\Text{Set}
\T(u)=(1+q^{2\La}\thi)\prod_{m=1}^nz_m(u/x_m-1)+
\sum_{m=1}^n{u\over x_m}\prod_{\tsize{k=1\atop k\ne m}}^n{u-x_k\over x_m-x_k}
\cdot(y_m^+ +\thi\>y_m^-)\,.
\endgather
$$
\Lm*
Coefficients of the \pol/ $\T(u)$ generate a commutative subalgera in
$\E(\Fo)$.
\endpro
\Lm{spectrq}
Let $\tau(u),\,Q(u)$ be a \sol/ to the \loc/. Set
$\Qo=\pi_1Q\lox\pi_nQ\ne 0$. Then $\T(u)\Qo=\tau(u)\Qo$.
Moreover, any \egv/ of $\T(u)$ has the form $\Qo$ for a suitable
\sol/ $\tau(u),\,Q(u)$ to the \loc/.
\endpro
Now let $q^2$ be a primitive \^{$N$-th} root of unity.
In this case we assume only that $q^{2\La_m}\ne 0$, $m=1\lc n$.
Define integers $d_1\lc d_n$ and $s_1\lc s_n$ as follows. If
$q^{4\La_m}\in\{\,1,\,q^2,\lc q^{2N-2}\}$ then define $d_m\in\{1\lc N-1\}$ by
$q^{4\La_m}=q^{2d_m}$ and set $s_m=d_m-1$. Otherwise set $d_m=N-1$, $s_m=N-1$.
We assume that $z_1\lc z_n$ are \wsep/, which means that all points
$z_mq^{2(s-\La_m)}\!$, $s=0\lc s_m$, and $z_mq^{2\La_m}$, $m=1\lc n$, are \pd/.
\par
Define sets $\Sc_1\lc \Sc_n$ as follows. If
$q^{4\La_m}\in\{\,1,\,q^2,\lc q^{2N-4}\}$ then define
$\Sc_m=\{\,z_mq^{2(s-\La_m)}\vert s=0\lc d_m\,\}$ similar to the case of
generic $q$. Otherwise define $\Sc_m=\{\,\ze_mq^{2s}\vert s=0\lc N-1\,\}$
with arbitrary $\ze_m\in\Co$. If $\ze_m\ne z_mq^{\pm2\La_m}$ then say that
$\Sc_m$ is {\it cyclic\/}. We assume that sets $\Sc_1\lc \Sc_n$ are \pdj/.
For instance, this is the case provided all of them are not cyclic.
\Lm{Q=Qr}
Let all sets $\Sc_1\lc \Sc_n$ be not cyclic or $\thi$ be generic.
Let $\tau(u)$ be a \sol/ to the \loc/. Then a \sol/ $Q\in\F$ to \eq/ \(Baxq)
is unique modulo a pseudoconstant factor.
\endpro
For any \pol/ $P(u)$ say that $P(u^N)$ is a quasiconstant.
\Lm{proporr}
Let $\tau(u)$ be a \sol/ to the \glob/. Then
\iitem{a)} $\ \deg\tau\le n$ and
$\tau(0)=(-1)^n(1+q^{2\La}\thi)\prod_{m=1}^nz_m$.
\iitem{b)} Let $\thi\nin\{\,q^{-2(s+\La)}\vert s=1\lc N-1\,\}$. Then for
a given $\tau(u)$, the required \pol/ \sol/ to \eq/ \(Baxq) is unique modulo
a quasiconstant factor.
\endpro
Let $\tau(u),\,Q(u)$ be a \sol/ to the \glob/ and $Q(u)$ has the smallest
possible degree. If $\tau(u)$ is also a \sol/ to the \loc/, then say that
$\tau(u)$ is an {\it \adm/} global \sol/.
\Th*
Let $q$ be a root of unity. Let $\thi$ be generic. Let $z_1\lc z_n$ be \wsep/.
Let $\Sc_1\lc\Sc_n$ be \pdj/. Then all claims of Theorem \[<>q] hold.
\endpro
\nt
The proof is completely similar to the case of generic $q$, because
all Lemmas \[Basolq]\,--\,\[spectrq] remain valid. A slight change is
necessary in the proof of Lemma \[Basolq] if some $\Sc_m$ is cyclic. Namely,
in this case $Q\vst{\Sc_m}\ne 0$ because \pol/s $Q(u)$ and $Q(q^2u)$ have no
common zeros. The last claim itself follows from the fact that $t_1\lc t_\ell$
is an \adm/ \sol/ to system \(Baeq).

\vsk>

\myRefs
\widest{KiR}

\ref\Key{B}
\by R\&J\&Baxter
\book Exactly solved models in statistical mechanics
\publ Academic Press \publaddr London
\yr 1982
\endref

\ref\Key{BIK}
\by N\&M\&Bogoliubov, A\&G\&Izergin and \Kor/
\book Quantum inverse scattering method and correlation functions
\publ Cambridge University Press
\yr 1993 \page 555
\endref

\ref\Key{CP1}
\by V\)\&Chari and A\&Pressley
\paper Yangians and \Rms/
\jour L'Enseignement Math. \vol 36 \yr 1990 \pages 267--302
\endref

\ref\Key{CP2}
\by V\)\&Chari and A\&Pressley
\paper Quantum affine allgebras
\jour \CMP/
\vol 142 \yr 1991 \pages 261--283
\endref

\ref\Key{DCK}
\by C\&De Concini and V\)\&G\&Kac
\paper Representations of quantum groups at roots of $1$
\inbook in Progress in Mathematics
\vol 92
\yr 1990 \pages 471--506
\publ Birkh\"auser
\endref

\ref \Key{EKS}
\by F\&Essler, \Kor/ and K\&Schoutens
\paper Fine structure of the \Ba/ for the spin-$1/2$ Heisenberg
{\sl XXX} model
\jour J. Phys. A \issue 15 \vol 25 \yr 1992 \pages 4115--4126
\endref

\ref\Key{F}
\by \Fadd/
\paper Lectures on quantum inverse scattering method
\inbook in Integrable Systems, Nankai Lectures on Math. Phys., 1987
\yr 1990 \pages 23--70
\ed X.-C\&Song \publ World Scientific \publaddr Singapore
\endref

\ref\Key{FT}
\by \Fadd/ and \Takh/
\paper
Quantum inverse problem method and the Heisenberg {\sl XYZ}-model
\jour Russian Math. Survey \vol 34 \yr 1979 \issue 5 \pages 11--68
\endref

\ref\Key{FT2}
\by \Fadd/ and \Takh/
\paper The spectrum and scattering of excitations in the one-dimen\-sional
isotropic Heisenberg model
\jour \ZNS/ \vol 109 \yr 1981 %\pages 134--178
\moreref
\jour \JSM/ \vol 24 \yr 1984 \pages 241--267
\endref

\ref\Key{Ki1}
\by \Kir/
\book Representation of quantum groups, combinatorics, $q$-orthogonal
polynomials and link invariants \bookinfo  Thesis
\yr 1990 \publ LOMI \publaddr Leningrad \page 300 \finalinfo (in Russian)
\endref

\ref\Key{Ki2}
\by \Kir/
\paper Combinatorical identities and the completness of states for
Heisenberg magnet
\jour \ZNS/ \vol 131 \yr 1983 \pages 88--105
\moreref
\jour \JSM/ \vol 30 \yr 1985 \pages 2298--3310
\endref

\ref
\by \refin\Kir/
\paper Completness of the states for the generalized Heisenberg model
\jour \ZNS/ \vol 134 \yr 1985 \pages 169--189
\moreref
\jour \JSM/ \vol 36 \yr 1987 \pages 115--128
\endref

\ref
\by \refin\Kir/ and N\&A\&Liskova
\paper Completeness of Bethe's states for generalized {\sl XXZ} model
\jour Preprint UTMS\~47 \yr 1994 \pages 1--20
\endref

\ref\Key{Ko}
\by \Kor/
\paper Calculation of norms of Bethe wave functions
\jour \CMP/ \vol 86 \yr 1982 \pages 391--418
\endref

\ref\Key{LS}
\by R\&P\)\&Langlands and Y\)\&Saint-Aubin
\paper Algebro-geometric aspects of the Bethe \eq/s
\jour Preprint \yr 1994 \pages 1--14
\endref

\ref\Key{NT}
\by \MN/ and \VT/
\paper Representations of Yangians with Gelfand-Zetlin bases
\jour Preprint \ifMag\else\adjust{\nl}\fi
UWS\~MRRS\~94\~148 \yr 1994 \pages 1--31
\endref

\ref\Key{RV}
\by \Resh/ and \Varch/
\paper Quasiclassical asymptotics of \sol/s to the \KZ/ \eq/s
\jour Preprint \yr 1994 \pages 1--29
\endref

\ref\Key{S1}
\by \Skl/
\paper Separation of \var/s in Gaudin model
\jour \ZNS/ \vol 164 \yr 1987 \pages 151--169
\moreref
\jour \JSM/ \vol 47 \yr 1989 \pages 2473--2488
\endref

\ref\Key{S2}
\by \Skl/
\paper Functional \Ba/
\inbook in Integrable and Superintegrable Systems
\ifMag\else\adjustnext{\hfill\nl}\fi
\yr 1990 \pages 8--33
\ed B\&A\&Kupershmidt
\publ World Scientific \publaddr Singapore
\endref

\ref\Key{S3}
\by \Skl/
\paper Quantum inverse scattering method. Selected topics
\inbook in Quantum Groups and Quantum Integrable Systems,
Nankai Lectures on Math. Phys., 1991
\yr 1992 \pages 63--97
\ed M.-L\&Ge \publ World Scientific \publaddr Singapore
\endref

\ref\Key{Sz}
\by G\&Szego
\book Orthogonal polynomials \yr 1939 \publ AMS
\endref

\ref\Key{T}
\by \VoT/
\paper Irreducible monodromy matrices for the \Rm/ of the {\sl XXZ}-model
and lattice local quantum Hamiltonians
\jour \TMP/ \vol 63 \yr 1985 \pages 440--454
\endref

\ref\Key{TV}
\by \VT/ and \Varch/
\paper Asymptotic \sol/s to the quantized \KZv/ \eq/ and \Bv/s
\jour Preprint UTMS 94\~46 \yr 1994 \pages 1--30
\endref

\endRefs

\bye